\documentclass [11pt]{article}
\usepackage[british]{babel}
\usepackage{amsfonts,amsmath,amssymb, enumerate,amsthm}
\usepackage{framed,enumitem} 
\usepackage{tikz}
\usepackage{graphics}
\usepackage[
bookmarksopen=true,
bookmarksopenlevel=1,
colorlinks=true,
linkcolor=darkblue,
linktoc=page,
citecolor=darkblue,
urlcolor=darkblue
]{hyperref}
\usepackage{mathdots}
\usepackage{xcolor}
\usepackage{diagbox}
\usepackage{colortbl}
\usepackage[absolute,overlay]{textpos}
\usepackage{graphicx}
\usetikzlibrary{calc}
\usetikzlibrary{decorations.pathreplacing}
\usetikzlibrary{positioning,patterns}
\usetikzlibrary{arrows,shapes,positioning}
\usetikzlibrary{decorations.markings}
\usepackage{mathtools}

\usepackage{bbm}
\definecolor{darkblue}{rgb}{0,0,0.5}

\newtheorem{theorem}{Theorem}
\newtheorem{prop}[theorem]{Proposition}
\newtheorem{lemma}[theorem]{Lemma}
\newtheorem{cor}[theorem]{Corollary}
\newtheorem{fact}[theorem]{Fact}
\newtheorem{claim}[theorem]{Claim}

\theoremstyle{definition}
\newtheorem{problem}[theorem]{Problem}
\newtheorem{conj}[theorem]{Conjecture}
\newtheorem{definition}[theorem]{Definition}
\newtheorem{remark}[theorem]{Remark}

\newcommand{\cE}{\mathcal{E}}
\newcommand{\cF}{\mathcal{F}}

\newcommand{\cI}{\mathcal{I}}
\newcommand{\cJ}{\mathcal{J}}

\newcommand{\cM}{\mathcal{M}}

\newcommand{\cP}{\mathcal{P}}
\newcommand{\cQ}{\mathcal{Q}}
\newcommand{\cR}{\mathcal{R}}
\newcommand{\cS}{\mathcal{S}}
\newcommand{\cT}{\mathcal{T}}
\newcommand{\cU}{\mathcal{U}}

\newcommand{\bN}{\mathbb{N}}

\newcommand{\bR}{\mathbb{R}}

\def\eps{{\varepsilon}}

\newcommand{\defn}{\emph}

\newcommand{\prob}[1]{\mathrm{\mathbb{P}}\left[#1\right]}
\newcommand{\cprob}[2]{\prob{#1 \;\middle|\; #2}}
\newcommand{\expn}[1]{\mathrm{\mathbb{E}}\left[#1\right]}

\def\lflr{\left\lfloor}
\def\rflr{\right\rfloor}
\def\lcl{\left\lceil}
\def\rcl{\right\rceil}

\newcommand{\IND}{\mathbbm{1}}

\newcommand{\op}{\operatorname}
\newcommand{\ff}[2]{\left\lfloor\frac{#1}{#2}\right\rfloor}
\newcommand{\cf}[2]{\left\lceil\frac{#1}{#2}\right\rceil}

\def\COMMENT#1{}
\def\TASK#1{}
\let\TASK=\footnote             
\let\COMMENT=\footnote          

\begin{document}

\title{On the factorization of Dirac hypergraphs}

\author{Yangyang Cheng\footnotemark[1]\thanks{Fakult\"at f\"ur Informatik und Mathematik, Universit\"at Passau, Innstra{\ss}e 41, 94032 Passau, Germany. \\ Research partially funded by the Deutsche Forschungsgemeinschaft (DFG, German Research Foundation) — 542321564. \\
\textit{Emails:} \href{mailto:yangyang.cheng@uni-passau.de}{\tt yangyang.cheng@uni-passau.de}, \href{mailto:stefan.glock@uni-passau.de}{\tt stefan.glock@uni-passau.de}, \href{mailto:andrew.lane@uni-passau.de}{\tt andrew.lane@uni-passau.de}.}
\and Stefan Glock\footnotemark[1]
\and Andrew Lane\footnotemark[1]}

\date{}

\maketitle

\begin{abstract}
A celebrated result from the 70s due to Baranyai states that the complete $k$-uniform hypergraph on $n$ vertices decomposes into perfect matchings if $k\mid n$. We prove a minimum degree version of this result. Note that a necessary condition for a $k$-uniform hypergraph $G$ on $n$ vertices to admit a decomposition into perfect matchings is that $G$ is vertex-regular and $k\mid n$. We show that, if this is satisfied and, in addition, $G$ has minimum $(k-1)$-degree at least $(\frac{1}{2}+o(1))n$, then it indeed decomposes into perfect matchings. The constant $1/2$ is optimal, as below this threshold, $G$ might not even contain a single perfect matching. This settles a conjecture of Glock, Kühn and Osthus from 2021. 
\end{abstract}

\section{Introduction}
A folklore construction shows that the edges of every complete graph of even order can be partitioned into $1$-factors, yielding a so-called $1$-factorization. An important class of structures in design theory, $1$-factorizations have significant theoretical importance and arise naturally in various scheduling problems, most notably in the organization of round-robin tournaments (see, e.g., \cite{CD:07}). In this paper, we consider $1$-factorizations of dense hypergraphs.

Given a positive integer $k$, a \defn{$k$-uniform hypergraph}, or \defn{$k$-graph} for short, is a hypergraph whose edges all have size $k$.
Given two hypergraphs $F$ and $G$, an \defn{$F$-decomposition of $G$} is a collection $\cF$ of subgraphs of $G$, each isomorphic to $F$, such that every edge of $G$ is contained in exactly one of the hypergraphs in $\cF$.
Graph and hypergraph decomposition problems involve determining necessary and sufficient conditions for a hypergraph $G$ to have an $F$-decomposition.
These problems have a long history in combinatorics: for example, Kirkman~\cite{kirkman:47} proved in 1847 that the complete graph $K_n$ has a $K_3$-decomposition if and only if $n \equiv 1,3 \pmod{6}$, and Walecki's theorem from 1892 states that for all odd $n$, $K_n$ has a decomposition into Hamilton cycles.
See, e.g., \cite{GKO:21} for a survey on graph and hypergraph decomposition problems.



A \defn{matching} in a hypergraph is a set of pairwise disjoint edges.
A \defn{perfect matching} or a \defn{$1$-factor} of a hypergraph is a matching that covers all vertices of the hypergraph. Baranyai's theorem, conjectured by Sylvester and proved by Baranyai in 1975~\cite{baranyai:75}, generalizes the folklore result on the existence of $1$-factorizations from graphs to hypergraphs. It states that for any integer $k\ge 2$ and any positive multiple $n$ of $k$, the complete $k$-uniform hypergraph $K_n^{(k)}$ admits a decomposition into perfect matchings, also called a \defn{$1$-factorization}, meaning that there exist pairwise edge-disjoint perfect matchings $M_1,\ldots,M_t$ of $K_n^{(k)}$ whose union is the entire edge set of $K_n^{(k)}$.

In the case of graphs ($2$-uniform hypergraphs), Baranyai's theorem has been generalized to graphs with high minimum degree.
Perkovic--Reed~\cite{PR:97} and H\"aggkvist (unpublished) proved that for any $\eps > 0$, there exists $n_0$ such that for any even $n \ge n_0$ and any $d \ge (\frac12 + \eps)n$, if $G$ is an $n$-vertex $d$-regular graph, then $G$ admits a $1$-factorization.
A celebrated result of Csaba, K\"uhn, Lo, Osthus, and Treglown~\cite{CBKLOT:16} determined the exact minimum degree threshold, resolving the long-standing $1$-factorization conjecture: they proved that for all sufficiently large even $n$ and $d \ge 2\lcl n/4\rcl-1$, every $n$-vertex $d$-regular graph has a $1$-factorization, and the bound on $d$ is best possible.

For hypergraphs, much less is known.
In fact, even guaranteeing just one perfect matching is already difficult.
K\"uhn and Osthus~\cite{KO:06} proved that for any $k\ge 2$ and $n$ divisible by $k$, any $k$-graph with minimum $(k-1)$-degree at least $(\frac12+o(1))n$ contains a perfect matching.
R\"odl, Ruci\'nski, and Szemer\'edi~\cite{RRS:09} later used the absorption method to determine the exact minimum value of $\delta_{k-1}(G)$ that ensures the presence of a perfect matching in a $k$-graph when $n$ is sufficiently large.
Research on this problem and its variants received much attention and went hand in hand with the systematic development of the absorption method (see, e.g.,~\cite{AGS:09,HPS:09,RRS:06,RRS:08,RRS:06b}).

It is then natural to try to generalize the above results for decompositions of graphs into perfect matchings to hypergraphs. Indeed, in their survey, Glock, K\"uhn, and Osthus~\cite[Conjecture~6.6]{GKO:21} conjectured a minimum degree version of Baranyai's theorem, stating that the minimum degree threshold guaranteeing one perfect matching is already sufficient for the existence of a $1$-factorization when also assuming the necessary condition of being vertex-regular.
Joos, K\"uhn, and Sch\"ulke~\cite{JKS:23} proved an approximate version of this conjecture, showing that an $r$-(vertex-)regular $k$-graph $G$ with $\delta_{k-1}(G)\ge (\frac12+\delta)n$ contains $(1-\delta)r$ edge-disjoint perfect matchings when $n$ is sufficiently large, as a corollary of their result on decompositions into cycle factors.

In this paper, we confirm the conjecture of Glock, K\"uhn, and Osthus~\cite{GKO:21}. To the best of our knowledge, this is the first result on hypergraph decompositions into spanning structures in the minimum degree setting.
Our proof employs the result of Joos, K\"uhn, and Sch\"ulke~\cite{JKS:23} as a key ingredient; specifically, we combine it with the absorption method to obtain the exact result.
%
Since the $k=2$ case is already known, we will only consider $k\ge 3$.
\begin{theorem} \label{thm:main}
Fix constants $k \ge 3$ and $\delta >0$. For all sufficiently large $n$ with $k \mid n$, if $G$ is a vertex-regular $k$-graph on $n$ vertices with $\delta_{k-1}(G) \ge (\frac12 + \delta)n$, then $G$ admits a $1$-factorization.
\end{theorem}
The constant $\frac12$ in Theorem~\ref{thm:main} is best possible.
Indeed, below this threshold, a $k$-graph satisfying $\delta_{k-1}(G) \ge (\frac12-o_k(1))n$ may not have a single perfect matching, and such a $k$-graph can be made vertex-regular and still not have a $1$-factorization (see Section~\ref{sec:lowerbound}).

\section{Preliminaries} \label{sec:preliminaries}
\subsection{Notation and partition definitions}
Let $\cP$ be a partition of a finite set $V$.
Say $\cP$ is \defn{$k$-equitable} if for all $U,U' \in \cP$, $k \mid |U|$ and $||U|-|U'||\in \{0,k\}$.
Equivalently, $\cP$ is $k$-equitable if each $U\in \cP$ has $|U| \in \left\{k\lflr \frac{|V|}{k|\cP|}\rflr,k\lcl\frac{|V|}{k|\cP|}\rcl\right\}$.
Say $\cP$ is \defn{equitable} if it is $1$-equitable.
For $\ell \in \bN$, an \defn{$\ell$-partition} is a partition of size $\ell$.
For each $W\subseteq V$, define $\cP[W] := \{U \in \cP:U \cap W \neq \emptyset\}$.
If $\cP$ and $\cQ$ are partitions of a set $V$, then $\cQ$ is a \defn{refinement} of $\cP$ if for all $Q \in \cQ$, there exists $P \in \cP$ such that $Q \subseteq P$.

For a graph $G$ and $\ell \in \bN$, an \defn{$\ell$-colouring} of $G$ is a partition of $V(G)$ into $\ell$ independent sets.
(We allow such a partition to be a multiset containing multiple copies of the empty set.)

For a $k$-graph $G$ and a collection $\cP$ of subsets of $V(G)$, let $G[\cP]$ be the subgraph of $G$ with vertex set $\bigcup_{U\in \cP}U$ and edge set $\bigcup_{U\in \cP}E(G[U])$.
In particular, if $\cP$ is a partition of $V(G)$ and $V\subseteq V(G)$, then $G[\cP[V]] = \bigcup_{\substack{U\in \cP \\ U \cap V \neq \emptyset}}G[U]$.

Given a hypergraph $G$ and sets $S,U \subseteq V(G)$, let $d_G(S) := |\{e\in E(G):S \subseteq e\}|$ be the \defn{degree} of $S$, let $d_{G,U}(S) := |\{e \in E(G):S\subseteq e \subseteq U\cup S\}|$, and denote $d_G(v) = d_G(\{v\})$ and $d_{G,U}(v) = d_{G,U}(\{v\})$ for $v\in V(G)$.
For $1\le \ell \le k$, let $\delta_{\ell}(G)$ and $\Delta_\ell(G)$ denote the minimum and maximum, respectively, of $d_G(S)$, taken over all $S \subseteq V(G)$ with $|S| = \ell$.
We refer to these as the \defn{minimum} and \defn{maximum $\ell$-degrees} of $G$, respectively.

Given variables $\eps, \delta>0$, we write $\eps \ll \delta$ to denote that given a fixed choice of $\delta > 0$, there exists $\eps_0 >0$ depending on $\delta$ such that for all $0 < \eps \le \eps_0$, the statement that follows holds.
When several variables are involved in a hierarchy, e.g.~$\mu \ll \eps \ll \delta$, the hierarchy is read right to left, so each variable depends on those to its right.
Say that a random event $E$ holds with high probability (w.h.p.) with respect to a variable $n$ if $\prob{E} \to 1$ as $n\to\infty$.

For convenience, in explanatory text, we refer to $k$-graphs $G$ on $n$ vertices that satisfy $\delta_{k-1}(G) \ge (\frac12+\delta) n$ for some constant $\delta$ as \defn{Dirac hypergraphs}.

\subsection{Required theorems}
We will use the Hajnal--Szemer\'edi theorem, which ensures the existence of equitable colourings with a number of colours close to its maximum degree.
\begin{theorem}[Hajnal--Szemer\'edi theorem~\cite{HS:70}] \label{thm:equitable}
    Let $G$ be a graph. For any integer $\ell > \Delta(G)$, $G$ has an equitable $\ell$-colouring.
\end{theorem}
The Hajnal--Szemer\'edi theorem can be easily generalized to $k$-equitable $\ell$-colourings.
\begin{cor} \label{cor:kequitable}
    Let $G$ be a graph, and let $n:=|V(G)|$. For all $k,\ell\in \bN$ with $k\mid n$ and $n,\ell > k\cdot \Delta(G)$, $G$ has a \defn{$k$-equitable $\ell$-colouring}.
\end{cor}
\begin{proof}
Since $n/k$ is an integer greater than $\Delta(G)$, by Theorem~\ref{thm:equitable}, there exists an equitable $n/k$-colouring $\Phi = \{V_1,\ldots,V_{n/k}\}$ of $G$.
Because $k \mid n$, we have $|V_i|=k$ for all $i \in [n/k]$.

Define the graph $G'$ by $V(G') = \Phi$ and $E(G') = \{\{V_i,V_j\}:E_G(V_i,V_j) \neq \emptyset\}$.
Observe that for $i \in [n/k]$, $d_{G'}(V_i) \le \sum_{v\in V_i}d_G(v) \le k\cdot \Delta(G)$.
Thus, $\ell>\Delta(G')$, so by Theorem~\ref{thm:equitable}, there exists an equitable $\ell$-colouring $\Phi' = \{U_1,\ldots,U_\ell\}$ of $G'$.
Let $\Phi''$ be the partition $\{W_1,\ldots, W_\ell\}$ of $V(G)$, where $W_i = \bigcup_{V\in U_i}V$ for all $i\in [\ell]$.
Then for $i,j \in [\ell]$, $|W_i| = k|U_i|$ is divisible by $k$, and since $\Phi'$ is equitable,
\[ ||W_i|-|W_j|| = |k|U_i|-k|U_j|| = k||U_i|-|U_j|| \in \{0,k\}. \]
Also, each $W_i$ is an independent set in $G$, since if $uv \in E(G[W_i])$, then either $uv \in E(V)$ for some $V\in \Phi$ (contradicting $\Phi$ being a colouring), or $uv \in E_G(V,V')$ for some distinct $V,V' \in \Phi$, in which case $\{V,V'\} \in E(G'[U_i])$ (contradicting $\Phi'$ being a colouring).
Therefore, $\Phi''$ is a $k$-equitable $\ell$-colouring of $G$, as desired.
\end{proof}

The following result of K\"uhn and Osthus shows that any Dirac $k$-graph has a perfect matching.
\begin{theorem}[K\"uhn, Osthus~\cite{KO:06}] \label{thm:rrsmatching}
    For all $k \ge 2$ and $\delta >0$, there exists $n_0$ such that for all $n \ge n_0$ with $k \mid n$, if $G$ is an $n$-vertex $k$-graph with $\delta_{k-1}(G) \ge (\frac12+\delta)n$, then $G$ contains a perfect matching.
\end{theorem}

We will also apply the following result of Glock, Gould, Joos, K\"uhn, and Osthus on uniformly random perfect matchings in Dirac $k$-graphs.
\begin{lemma}[Glock, Gould, Joos, K\"uhn, Osthus~{\cite[Lemma 4.1]{GGJKO:21}}] \label{lem:uniformmatching}
    Let $1/n \ll 1/C \ll \delta,1/k$, where $k\ge 2$ and $k\mid n$.
    Let $G$ be an $n$-vertex $k$-graph with $\delta_{k-1}(G) \ge (\frac12+\delta)n$, and let $M$ be a uniformly random perfect matching of $G$.
    Then, for each $e\in E(G)$, we have
    \[ \frac{1}{Cn^{k-1}} \le \prob{e\in E(M)} \le \frac{C}{n^{k-1}}. \]
\end{lemma}

Joos, K\"uhn, and Sch\"ulke proved that an almost-regular Dirac $k$-graph has an approximate $1$-factorization.
Say a $k$-graph $G$ is \defn{$\varrho$-almost $r$-regular} if $d_G(v) = (1\pm \varrho)r$ for all $v\in V(G)$.
\begin{theorem}[Joos, K\"uhn, Sch\"ulke~\cite{JKS:23}] \label{thm:approxdecomp}
    For all $k\ge 2$ and $\delta>0$, there exist $n_0, \varrho>0$ such that every $k$-graph $G$ on $n\ge n_0$ vertices with $k \mid n$ and $\delta_{k-1}(G) \ge (\frac12 + \delta)n$ which is $\varrho$-almost $r$-regular for some $r\ge 0$ contains $(1-\delta)r$ pairwise edge-disjoint perfect matchings.
\end{theorem}

We will use the following form of the Chernoff bound.
\begin{theorem}[See, e.g.,~{\cite[Theorem 2.10]{JLR:00}}] \label{thm:chernoffhyper}
    Let $X$ be a finite sum of independent hypergeometric or Bernoulli random variables, and let $\mu := \expn{X}$.
    Then for any $t>0$, we have $\prob{X \ge \mu+t} \le e^{-t^2/2(\mu+t/3)}$ and $\prob{X \le \mu-t} \le e^{-t^2/2\mu}$.
\end{theorem}

\section{Proof outline and main lemmas} \label{sec:outline}
We prove Theorem~\ref{thm:main} using a variant of the iterative absorption method, in the form of iteratively refined vertex partitions as developed by Barber, K\"uhn, Lo, and Osthus~\cite{BKLO:16}.
One of our main goals is to set aside an absorbing structure $A$ within the $k$-graph $G$ with the property that any sufficiently small leftover subgraph of $G$ admits a $1$-factorization when combined with $A$.
In our case, we require this small leftover to be contained in $G[\cU] = \bigcup_{U\in \cU}G[U]$ for some partition $\cU$ of $V(G)$ into constant-size sets, meaning that each edge of the leftover is contained entirely in one of the constant-size vertex sets in $\cU$.
The following lemma details the key properties of the absorbing structure $A$.
(Here, $K_n^{(k)}$ is assumed to have the same vertex set as $G$.)
\begin{lemma}[Absorption lemma] \label{lem:absorbing}
Let $k\ge 3$.
Let $1/n \ll 1/m \ll \gamma \ll \delta$ such that $k\mid m$ and $k \mid n$.
Let $G$ be an $n$-vertex $k$-graph with $\delta_{k-1}(G) \ge (\frac12+\delta)n$.
Let $\cU$ be a partition of $V(G)$ into $m$- and $(m+k)$-sets, and let $d_0 \in \bN$.
Then there exists a spanning subgraph $A$ of $G-G[\cU]$ such that $\Delta_{k-1}(A) \le \gamma n$ and for any $d_0$-regular spanning subgraph $K$ of $K_n^{(k)}[\cU]$, $A \cup K$ admits a $1$-factorization.
\end{lemma}
We prove Lemma~\ref{lem:absorbing} through a series of reductions, which involves a novel concept that we call \defn{compatible matchings}.
The first and quite straightforward reduction, which appears in the proof of Lemma~\ref{lem:absorbing} itself, is that it suffices to construct an absorber $A_U$ for each $U \in \cU$ individually, where $|V(A_U)| = o(n)$ and $A_U \cup H_U$ has a $1$-factorization for any $d_0$-regular $H_U \subseteq G[U]$, because we can then find such edge-disjoint absorbers one by one and ``fill in" with matchings on $V(G)\setminus V(A_U)$.
The second reduction, given by Lemma~\ref{lem:absorbingperfect}, is that it suffices to construct $A_U$ such that for all $H_U \subseteq G[U]$, $A_U \cup H_U$ has a proper edge-colouring (i.e., a decomposition into matchings) in which, given a fixed partition $\cP$ of $V(A_U)$ into $k$-sets, any two vertices in the same set in $\cP$ are incident to the same colours; this is weaker than the notion of a $1$-factorization, since colour classes need not span the entire vertex set.
We call the colour classes of such an edge-colouring \defn{$\cP$-compatible matchings}.
We accomplish this reduction by expanding the vertex set of the absorber $A_U$ and appending an absorbing gadget to each $P\in \cP$ using these new vertices, so that each $\cP$-compatible matching is completed to a perfect matching.

The third reduction, given by Lemma~\ref{lem:absorbingcompatible}, is that it suffices to construct $A_U$ using ``fake edges", i.e., $A_U$ can be any $k$-(multi)graph on an appropriate vertex set, not necessarily a subgraph of $G$.
This follows because each ``fake edge" can be replaced by actual edges in $G$ by expanding the vertex set of $A_U$ and adding a matching on the fake edge along with some new vertices that functions exactly like the fake edge with respect to $\cP$-compatible matchings.
Finally, in Lemma~\ref{lem:absorbingmulti}, we construct $A_U$ using fake edges by adding $\cP$-compatible matchings that can absorb any $d_0$-regular subgraph of $G[U]$.
We present these reductions in the opposite order, instead describing a multi-step process in which an arbitrary $d_0$-regular subgraph of $G[U]$ is absorbed by fixed absorbers into structures with increasingly strong properties (and increasingly large vertex sets), eventually yielding a $1$-factorization.

While Theorem~\ref{thm:approxdecomp} allows us to decompose almost all edges of $G-A$ into perfect matchings, it does not guarantee that the leftover edges are contained in $G[\cU]$ for a partition $\cU$ of $V(G)$ into small sets.
We will thus use the iterative absorption method to obtain an approximate decomposition of $G-A$ into perfect matchings so that the leftover is forced into smaller and smaller vertex sets in successively refined partitions.
To this end, we have the following definition of a vertex partition in which the sets in the partition approximately preserve the minimum $(k-1)$-degree and vertex-regularity conditions of $G$.
\begin{definition}
    Given a $k$-graph $G$, say that a partition $\cP$ of $V(G)$ is a \defn{$(\delta,\beta,L,\lambda)$-partition} for $G$ if the following hold for all $U \in \cP$:
    \begin{enumerate}[label=(P\arabic*)]
        \item $\cP$ is a $k$-equitable $L$-partition of $V(G)$;
        \item For all $S\in \binom{V(G)}{k-1}$, $d_{G,U}(S) \ge (\frac12+\delta -  |U|^{-1/3})|U|$;
        \item $G[U]$ is $\lambda|U|^{-1/3}$-almost $\beta |U|^{k-1}$-regular.
    \end{enumerate}
\end{definition}
The following lemma asserts the existence of $(\delta,\beta,L,\lambda)$-partitions of almost-regular Dirac $k$-graphs, which we will later apply inductively to show the existence of successively refined partitions with these same properties.
\begin{lemma} \label{lem:vortexpartition}
    Let $k \ge 3$ and $1/n \ll 1/L\ll \delta<1$ with $k \mid n$.
    Suppose $G$ is an $n^{-1/3}$-almost $r$-regular $k$-graph on $n$ vertices with $\delta_{k-1}(G) \ge (\frac12+\delta-n^{-1/3})n$.
    Then there exists a $(\delta, r/n^{k-1},L,1)$-partition $\cP$ for $G$.
\end{lemma}
The following cover-down lemma shows that given such a partition $\cP$, we can decompose all edges in $G-G[\cP]$ into perfect matchings by borrowing a small number of edges from $G[\cP]$.
To achieve this, we apply Theorem~\ref{thm:approxdecomp} to decompose almost all edges of $G-G[\cP]$ into perfect matchings and then apply Corollary~\ref{cor:kequitable} to an auxiliary graph in order to separate the remaining edges in $G-G[\cP]$ into matchings, each of which can be completed to a perfect matching using randomly chosen matchings (see Section~\ref{sec:sparsematchings}). A more detailed proof sketch can be found in Section~\ref{sec:coverdown}.

\begin{lemma}[Cover-down lemma] \label{lem:coverdown}
    Let $k \ge 3$ and $1/n \ll \eta \ll 1/L \ll \delta < 1$ with $k\mid n$.
    Suppose that for some $r$ and $\beta$, $G$ is an $n^{-1/3}$-almost $r$-regular $k$-graph on $n$ vertices with $\delta_{k-1}(G) \ge (\frac12+\delta)n$, $\cP$ is a $(\delta,\beta,L,2)$-partition for $G$, and for all $U \in \cP$, we have $k \mid \sum_{v\in U}d_G(v)$.
    Then there exists a subgraph $H$ of $G$ such that $G-G[\cP] \subseteq H$, $\Delta_{k-1}(H[\cP]) \le \eta n$, and $H$ admits a $1$-factorization.
\end{lemma}

Now, after applying Lemma~\ref{lem:coverdown}, the leftover edges, i.e., the edges that have not been decomposed into perfect matchings, have been forced into a subgraph of $G[\cP]$.
Further, for each $U\in \cP$, the leftover edges in $G[U]$ approximately inherit the degree conditions from $G$.
Thus, assuming $G[U]$ has its own partition $\cP_U$ constructed using Lemma~\ref{lem:vortexpartition}, we can apply Lemma~\ref{lem:coverdown} to $G[U]$ again to obtain more perfect matchings of $G$ and force the leftover edges into even smaller partition sets.
(We need to be careful that we did not borrow any edges from $G[\cP_U]$ in an earlier step, since otherwise, we could gradually destroy the minimum $(k-1)$-degree condition as our partition sets become smaller, but we can do this by simply deleting $G[\cP_U]$ from $G$ earlier on and adding it back when appropriate.)
Continuing in this fashion, we will iteratively absorb almost all edges into perfect matchings until the leftover edges are contained in $G[\cU]$ for a partition $\cU$ of $V(G)$ into constant-size sets.

Putting this all together, we construct a $1$-factorization of $G$ in the following manner.
First, we apply Lemma~\ref{lem:vortexpartition} iteratively to obtain a sequence of successively refined partitions $\cP_0,\cP_1,\ldots,\cP_t$, in which all vertex sets approximately inherit the minimum $(k-1)$-degree and vertex-regularity of $G$.
We then apply Lemma~\ref{lem:absorbing} to set aside an absorber $A$ which can absorb any leftover in $G[\cP_t]$.
After that, we apply Lemma~\ref{lem:coverdown} iteratively to decompose all edges outside of $G[\cP_t]$ into perfect matchings.
Our absorber $A$, combined with the leftover in $G[\cP_t]$, admits a $1$-factorization, and so we have decomposed all edges of $G$ into perfect matchings.

\subsection{Layout and organization}
The remainder of the paper is organized as follows.
In Section~\ref{sec:sparsematchings}, we prove some auxiliary lemmas on randomly chosen pairwise edge-disjoint matchings on a fixed vertex set.
In Section~\ref{sec:vortex}, we prove Lemma~\ref{lem:vortexpartition} and apply it to obtain a sequence of iteratively refined partitions in Lemma~\ref{lem:vortex}.
In Section~\ref{sec:coverdown}, we prove the Cover-down lemma (Lemma~\ref{lem:coverdown}) and apply it to the sequence of partitions from Section~\ref{sec:vortex} in Lemma~\ref{lem:iterative}.
In Section~\ref{sec:absorption}, we prove the Absorption lemma (Lemma~\ref{lem:absorbing}).
In Section~\ref{sec:main}, we combine these intermediate results to prove Theorem~\ref{thm:main}.
In Section~\ref{sec:lowerbound}, we show the existence of vertex-regular $k$-graphs with minimum $(k-1)$-degree at least $(\frac12-o(1))n$ but no $1$-factorization. In Section~\ref{sec:conclusion}, we discuss potential avenues for future research.

\section{Matchings with small collective $(k-1)$-degree} \label{sec:sparsematchings}
The following lemma shows that if we choose a collection of edge-disjoint matchings randomly, each covering all but a small constant proportion of the vertices of a Dirac $k$-graph and collectively making up a small constant proportion of the edges, then the union of these matchings has small $(k-1)$-degree.
We will use it to complete certain collections of partial matchings to perfect matchings, maintaining edge-disjointness.

\begin{lemma} \label{lem:sparsematchings}
    Let $k\ge 3$, let $1/n \ll \delta_1 \ll \delta_2 \ll \delta$ with $k \mid n$, and let $\ell \le \delta_1 n^{k-1}$.
    Suppose $G$ is an $n$-vertex $k$-graph with $\delta_{k-1}(G) \ge (\frac12 + \delta)n$ and $V_1,\ldots, V_\ell$ are subsets of $V(G)$ such that $|V_i| \le \delta_1n$ and $k$ divides $n-|V_i|$ for all $i\in [\ell]$.
    Choose $M_1,\ldots,M_\ell$ randomly and iteratively so that for each $i \in [\ell]$, letting $F_{i-1} := \bigcup_{j=1}^{i-1}M_j$, $M_i$ is a uniformly random perfect matching of $G - F_{i-1} - V_i$ if $\Delta_{k-1}(F_{i-1}) \le \delta_2n$ and is empty otherwise.
    Then w.h.p. as $n\to\infty$, $M_1,\ldots,M_\ell$ are edge-disjoint matchings of $G$ with $\Delta_{k-1}(\bigcup_{i=1}^\ell M_i) \le \delta_2 n$ and $V(M_i) = V(G)\setminus V_i$ for all $i\in [\ell]$.
\end{lemma}

\begin{proof}[Proof of Lemma~\ref{lem:sparsematchings}]
Let $F_\ell := \bigcup_{i=1}^{\ell}M_i$.
Let $S \in \binom{V(G)}{k-1}$.
For $i \in [\ell]$, let $X_{S,i} := d_{M_i}(S)$ and $Y_{S,i} := d_{F_i}(S)= \sum_{j=1}^i X_{S,i}$, and let $Y_{S,0} := 0$.
By Lemma~\ref{lem:uniformmatching} with $C=(1-\delta_1)^{k-1}/\delta_2 \gg 1/\delta$, since $G-F_{i-1}-V_i$ and $G-V_i$ have order at least $(1-\delta_1)n$ for all $i$,
\begin{equation} \label{eq:sparsematchings1}
    \cprob{X_{S,i}=1}{M_1,\ldots,M_{i-1}} = \cprob{\bigcup_{\substack{e\in E(G) \\ S\subseteq e}}\{e\in E(M_i)\}}{M_1,\ldots,M_{i-1}} \le \frac{d_G(S)}{\delta_2 n^{k-1}} \le \frac{1}{\delta_2 n^{k-2}}.
\end{equation}
Let $p := \frac{1}{\delta_2 n^{k-2}}$, let $Z_1,\ldots,Z_{\lflr \delta_1 n^{k-1} \rflr}$ be independent Bernoulli random variables with $\prob{Z_i=1}=p$ for all $i$, and for $i\in [0,\delta_1n^{k-1}]$, let $B_i := \sum_{j=1}^i Z_i \sim \op{Bin}(i,p)$.
It follows from~\eqref{eq:sparsematchings1} that $B_i$ has (first order) stochastic dominance over $Y_{S,i}$ for all $i \in [0,\ell]$, i.e., for all $x\in \bN_0$, $\prob{Y_{S,i}\ge x} \le \prob{B_i \ge x}$.

Letting $B := B_{\lflr \delta_1 n^{k-1}\rflr}$, we have $\expn{B} = (1+o(1))\frac{\delta_1 n^{k-1}}{\delta_2 n^{k-2}} \le \frac{\delta_2n}{2}$ and $B_i \le B$ for all $i \in [0,\ell]$, so by Theorem~\ref{thm:chernoffhyper}, since $1/n \ll \delta_1 \ll \delta_2$,
\[
    \prob{\Delta_{k-1}(F_\ell) > \delta_2n} \le \sum_{S\in \binom{V(G)}{k-1}}\prob{d_{F_\ell}(S) > \delta_2n} 
    \le \binom{n}{k-1}\prob{B > \delta_2 n} 
    \le n^{k-1}e^{-\Omega\left(\delta_2 n\right)} = o(1).
\]
Thus, with high probability, $\Delta_{k-1}(F_\ell) \le \delta_2 n$.
In particular, $\Delta_{k-1}(F_i) \le \delta_2 n$ for all $i \in [0,\ell]$, so by the definition of $M_1,\ldots,M_\ell$, this implies that each $M_i$ is a perfect matching of $G-F_{i-1}-V_i$.
Therefore, $M_1,\ldots,M_\ell$ are pairwise edge-disjoint matchings in $G$ such that $V(M_i) = V(G) \setminus V_i$ for all $i \in [\ell]$ and $\Delta_{k-1}(\bigcup_{i=1}^\ell M_i) \le \delta_2 n$, as desired.
\end{proof}

We will also need an extension of Lemma~\ref{lem:sparsematchings} that shows the existence of matchings that are ``uniform" with respect to a fixed partition of the vertices of a $k$-graph, in the sense that there are not too few edges with a fixed intersection size with each set in the partition.
\begin{definition} \label{def:types}
For $L,k \in \bN$, let $\cI_{L,k} := \{I \in \bN_0^L:\sum_{i=1}^L I_i = k\}$ and $\cJ_{L,k} := \{I \in \cI_{L,k}:\max_i I_i < k\}$.
Observe that $|\cI|,|\cJ|$ depend only on $L$ and $k$; in particular, the function $[L]^{[k]}\to \cI_{L,k}$ that maps $f\in [L]^{[k]}$ to $(|f^{-1}(1)|,\ldots,|f^{-1}(L)|)$ is surjective, so $|\cJ_{L,k}| \le |\cI_{L,k}| \le L^k$.
\end{definition}

\begin{lemma} \label{lem:uniformtypes}
    Let $k \ge 3$, let $1/n \ll \eps_0 \ll \eps_1 \ll \eps_2 \ll 1/L \ll \delta$ with $k \mid n$, and let $d := \lflr \eps_1 n^{k-1}\rflr$ and $\alpha := 1/L$.
    Suppose $G$ is an $n$-vertex $k$-graph, and let $\cP=\{U_1,\ldots,U_L\}$ be a $k$-equitable partition of $V(G)$ such that $d_{G,U}(S) \ge (\frac12+\delta)|U|$ for all $U \in \cP$ and $S \in \binom{V(G)}{k-1}$.
    Then there exists a $d$-regular spanning subgraph $F$ of $G-G[\cP]$ such that $\Delta_{k-1}(F) \le \eps_2n$, and for all $J \in \cJ_{L,k}$, there are at least $\eps_0 n^k$ edges $e$ in $F$ with $|e\cap U_i|=J_i$ for all $i \in [L]$.
\end{lemma}
\begin{proof}
Let $\eps$ be a new constant with $\eps_0 \ll \eps \ll \eps_1$, and let $G' := G-G[\cP]$.
Observe that $\delta_{k-1}(G) \ge (\frac12 +\delta)n$ and $\Delta_{k-1}(G[\cP]) \le \max_{U\in \cP}|U| \le 2\alpha n \le \delta n /2$, so $\delta_{k-1}(G') \ge (\frac12+\frac\delta2)n$.
Choose $M_1,\ldots,M_\ell$ randomly as in Lemma~\ref{lem:sparsematchings} with $\ell=d$ and $V_1=\cdots=V_\ell = \emptyset$, and let $F := \bigcup_{i=1}^d M_i$.
By Lemma~\ref{lem:sparsematchings} applied to $G'$ with $\delta$ replaced by $\delta/2$, w.h.p. $\Delta_{k-1}(F) \le \eps_2 n$.

Let $\cI = \cI_{L,k}$ and $\cJ := \cJ_{L,k}$ as in Definition~\ref{def:types}.
Define $g:2^{V(G)}\to \bN_0^L$ by $g(S) = (|S\cap U_1|,\ldots,|S\cap U_L|)$, and for $J \in \cJ$ and $K \in \{G,G',F\}$, let $K_J$ be the spanning subgraph of $K$ with $E(K_J) = E(K) \cap g^{-1}(J)$.
We are done if we can show that with high probability, $|E(F_J)| \ge \eps_0 n^k$ for all $J \in \cJ$.

Let $J\in \cJ$, and note that $G'_J = G_J$.
Let $j \in [L]$ with $J_j >0$, and let $I \in \bN_0^L$ with $I_j = J_j-1$ and $I_i = J_i$ for all $i \neq j$.
Let $S \in g^{-1}(I)$, and note that 
\[d_{G'_J}(S) = d_{G,U_j}(S) \ge \left(\frac12+\delta\right)|U_j| \ge \left(\frac12+\frac\delta2\right)\alpha n.\]
For $i\in [d]$, let $X_{S,i} := d_{M_i,U_j}(S)+\IND_{M_i = \emptyset}$.
Observe that the events $\{S \cup \{v\} \in E(M_i)\}$ are pairwise disjoint for all $v\in V(G)$.
If $\Delta_{k-1}(F_{i-1}) > \eps_2n$, then $M_i$ is empty, so $X_{S,i}=1$, and otherwise, by Lemma~\ref{lem:uniformmatching} with $C = 1/2\eps_2 \gg 1/\delta$,
\begin{align*}
    \cprob{X_{S,i}=1}{M_1,\ldots,M_{i-1}} &= \sum_{v\in U_j}\cprob{S \cup \{v\} \in E(M_i)}{M_1,\ldots,M_{i-1}} \\
    &\ge \frac{d_{G-F_{i-1},U_j}(S)}{Cn^{k-1}} \ge \frac{2\eps_2}{n^{k-1}}\left(\frac12 + \frac\delta2 - \eps_2\alpha^{-1}\right)\alpha n \ge \frac{\eps_2 \alpha}{n^{k-2}}.
\end{align*}
Letting $B \sim \op{Bin}(d,\frac{\eps_2\alpha}{n^{k-2}})$ and $Y_S := \sum_{i=1}^{d}X_{S,i}$, it follows that $Y_S$ has (first order) stochastic dominance over $B$.
Also, $\expn{B} = (1+o(1))\eps_1\eps_2 \alpha n \ge 2 \eps n$.
Thus, by Theorem~\ref{thm:chernoffhyper},
\[
    \prob{\bigcup_{S\in g^{-1}(I)}\{Y_S < \eps n\}} \le \sum_{S \in g^{-1}(I)}\prob{B < \eps n} \le n^{k-1}e^{-\Omega(\eps n)}.
\]
Also, note that if $\Delta_{k-1}(F) \le \eps_2 n$ and $Y_S \ge \eps n$ for all $S \in g^{-1}(I)$, then all $M_i$ are non-empty, so $Y_S = d_{F,U_j}(S)$ for all $S$, and so by counting pairs $(S,v) \in g^{-1}(I) \times U_j$ such that $S \cup \{v\} \in E(F)$, we have
\[ k|E(F_J)| \ge \sum_{S \in g^{-1}(I)}d_{F,U_j}(S) = \Omega_k\left((\alpha n)^{k-1}\eps n\right) \ge k \eps_0n^k, \]
implying that $|E(F_J)| \ge \eps_0 n^k$.
This means that $\prob{|E(F_J)| < \eps_0 n^k,\, \Delta_{k-1}(F) \le \eps_2 n} \le n^{k-1}e^{-\Omega(\eps n)}$.
Further, as noted in Definition~\ref{def:types}, $|\cJ| \le L^k$.
Thus,
\begin{align*}
    \prob{\bigcup_{J \in \cJ}\left\{|E(F_J)| < \eps_0 n^k\right\}} &\le \prob{\bigcup_{J\in \cJ}\left\{|E(F_J)| < \eps_0 n^k\right\},\Delta_{k-1}(F) \le \eps_2 n}+\prob{\Delta_{k-1}(F)>\eps_2n} \\
    &\le L^kn^{k-1}e^{-\Omega(\eps n )}+o(1) = o(1),
\end{align*} 
so w.h.p. $\Delta_{k-1}(F) \le \eps_2 n$ and $|E(F_J)| \ge \eps_0 n^k$ for all $J\in \cJ$, as desired.
\end{proof}

\section{Vortex lemma} \label{sec:vortex}
In this section, we prove Lemma~\ref{lem:vortexpartition} and apply it iteratively to obtain a sequence of successively refined partitions of the vertices of a vertex-regular Dirac $k$-graph so that each partition set induces an almost-regular Dirac $k$-graph.
\subsection{Proof of Lemma~\ref{lem:vortexpartition}}
We begin by proving Lemma~\ref{lem:vortexpartition}.
We prove the existence of a partition of the vertex set with the required vertex-degree and $(k-1)$-degree conditions by choosing a uniformly random $k$-equitable $L$-partition of the vertex set and showing that with high probability, these degrees are all close to their expected values.
We will apply the following concentration inequality in order to do so.

For $c>0$ and $n\in \bN$, say that a function $f:\{0,1\}^n\to \bR$ is \defn{$c$-Lipschitz} if changing $x\in \{0,1\}^n$ on one index changes $f(x)$ by at most $c$, i.e., for all $i \in [n]$ and $x\in \{0,1\}^n$,
\[ |f(x_1,\ldots,x_{i-1},1,x_{i+1},\ldots,x_n)-f(x_1,\ldots,x_{i-1},0,x_{i+1},\ldots,x_n)| \le c. \]
\begin{theorem}[McDiarmid's inequality~\cite{mcdiarmid:89}] \label{thm:mcdiarmid}
    Let $X_1,\ldots,X_n$ be independent random variables taking values in $\{0,1\}$, and let $c,t>0$.
    Then for any $c$-Lipschitz function $f:\{0,1\}^n\to \bR$,
    \[ \prob{f(X_1,\ldots,X_n) \ge \expn{f(X_1,\ldots,X_n)}+t} \le e^{-t^2/2c^2n} \]
    and
    \[ \prob{f(X_1,\ldots,X_n) \le \expn{f(X_1,\ldots,X_n)}-t} \le e^{-t^2/2c^2n}. \]
\end{theorem}
\begin{proof}[Proof of Lemma~\ref{lem:vortexpartition}]
Let $\beta := r/n^{k-1}$ and $V:=V(G)$.
Note that $\beta =\Theta_k(1)$ by the minimum $(k-1)$-degree condition on $G$.
Let $\cP = \{U_1,\ldots,U_L\}$ be a random $k$-equitable partition of $V$ chosen such that each $U_i$ is a uniformly random subset of $V$ of a fixed size (which exists since $k \mid n$).
We will show that with non-zero probability, $\cP$ is a $(\delta,\beta,L,1)$-partition for $G$.

Fix a set $U \in \cP$, and let $m := |U|$.
Let $U'$ be a random subset of $V$ such that each $v\in V$ is in $U'$ with probability $p:=m/n = L^{-1}\pm kn^{-1}$, independently of all other vertices.
Observe that, when conditioned on $|U'|=|U|$, $U$ and $U'$ are equal in distribution.
Also, by standard binomial coefficient estimates,
\begin{equation} \label{eq:vortexbinprob}
    \prob{|U'|=m} = \binom{n}{m}p^{m}(1-p)^{n-m} =\Omega\left(\sqrt{\frac{n}{m(n-m)}}\right) = \Omega(m^{-1/2}).
\end{equation}
For each $A\subseteq V$ with $|A| \le k-1$, define $f_A:\{0,1\}^V\to \bR$ by
\[ f_A\left((x_u)_{u\in V}\right) = \sum_{\substack{e\in E(G) \\ A\subseteq e}}\prod_{u \in e\setminus A}x_u. \]
Observe that for $x\in \{0,1\}^V$, if $X= \{u \in V:x_u=1\}$, then $f_A(x) = d_{G,X}(A)$.
It follows that changing $x$ on a single index $u$ changes $f_A(x)$ by at most $d_{G}(A\cup \{u\})$, so $f_A$ is $c$-Lipschitz for $c= \Delta_{|A|+1}(G)$.
For $u\in V$, let $X_u := \IND_{u\in U'}$, and let
\[ Y_A := f_A((X_u)_{u\in V}) = d_{G,U'}(A). \]
For $A \subseteq V$ of size $\ell \le k$,
\begin{equation} \label{eq:vortexdegprecise}
    \expn{Y_A} = \sum_{\substack{e\in E(G) \\ A\subseteq e}}\expn{\prod_{u\in e\setminus A}\IND_{u\in U'}} = p^{k-\ell}d_{G}(A).
\end{equation}
Let $W \in \{U,U'\}$; for $S\in \binom{V}{k-1}$, let $A_{W,S}$ be the event that $d_{G,W}(S) < (\frac12+\delta - m^{-1/3})m$; and for $v\in V$, let $A_{W,v}$ be the event that $|d_{G,W}(v) - \beta m^{k-1}|> \beta m^{k-4/3}$.
If $U$ fails to satisfy (P2)--(P3), then $A_{U,S}$ or $A_{U,v}$ holds for some $S \in \binom{V}{k-1}$ or $v \in V$.

\begin{claim} \label{clm:vortexk-1set}
    We have $\sum_{S \in \binom{V}{k-1}} \prob{A_{U,S}} < \frac{1}{2L}$.
\end{claim}
\begin{proof}
Let $S\in \binom{V}{k-1}$.
By \eqref{eq:vortexdegprecise}, since $\delta_{k-1}(G) \ge (\frac12+\delta - n^{-1/3})n$,
\begin{align*}
    \expn{Y_S} \ge \left(\frac12+\delta - n^{-1/3}\right)m,
\end{align*}
so
\[ \expn{Y_S}-\left(\frac12+\delta - m^{-1/3}\right)m \ge \left(1-p^{1/3}\right)m^{2/3}. \]
Because $f_S$ is $1$-Lipschitz (since $\Delta_k(G) \le 1$), it follows from Theorem~\ref{thm:mcdiarmid} that
\[  \prob{A_{U',S}} \le \prob{Y_S < \expn{Y_S}-(1-p^{1/3})m^{2/3}} \le e^{-\Omega\left(\frac{m^{4/3}}{m}\right)}. \]
Also, by \eqref{eq:vortexbinprob},
\[ \prob{A_{U,S}} = \cprob{A_{U',S}}{|U'|=|U|} \le \frac{\prob{A_{U',S}}}{\prob{|U'|=|U|}} \le m^{1/2}e^{-\Omega\left(\frac{m^{4/3}}{m}\right)}. \]
Therefore,
\[
        L\sum_{S \in \binom{V}{k-1}}\prob{A_{U,S}} \le Ln^{k-1}m^{1/2} e^{-\Omega\left(m^{1/3}\right)} <\frac12,
\]
completing the proof of Claim~\ref{clm:vortexk-1set}.
\end{proof}

\begin{claim} \label{clm:vortexvertex}
    We have $\sum_{v\in V}\prob{A_{U,v}} < \frac{1}{2L}$.
\end{claim}
\begin{proof}[Proof of Claim]
Let $v\in V$, $Y_v := Y_{\{v\}} = d_{G,U'}(v)$, and $f_v := f_{\{v\}}$.
By \eqref{eq:vortexdegprecise}, since $G$ is $n^{-1/3}$-almost $\beta n^{k-1}$-regular,
\[
    \expn{Y_v} = p^{k-1}\left(1\pm n^{-1/3}\right)\beta n^{k-1}= \left(1\pm n^{-1/3}\right)\beta m^{k-1}.
\]
Thus, because $f_v$ is a $c$-Lipschitz function for $c= \Delta_2(G) = O(n^{k-2})$, it follows from Theorem~\ref{thm:mcdiarmid} that
\[ \prob{A_{U',v}} \le \prob{|Y_v-\expn{Y_v}| > (1-p^{1/3})\beta m^{k-4/3}} \le \exp\left(-\Omega\left(\frac{m^{2(k-4/3)}}{n^{2(k-2)}n}\right)\right). \]
Also, by \eqref{eq:vortexbinprob},
\[ \prob{A_{U,v}} = \cprob{A_{U',v}}{|U'|=|U|} \le \frac{\prob{A_{U',v}}}{\prob{|U'|=|U|}} \le m^{1/2}e^{-\Omega(p^{2k-3}m^{1/3})}. \]
Therefore,
\[ L\sum_{v\in V}\prob{A_{U,v}} \le L n m^{1/2}e^{-\Omega(p^{2k-3} m^{1/3})} < \frac12, \]
completing the proof of Claim~\ref{clm:vortexvertex}.
\end{proof}
By Claims~\ref{clm:vortexk-1set}~and~\ref{clm:vortexvertex}, the probability that some $U \in \cP$ fails to satisfy (P2)--(P3) is at most
\begin{align*}
    &\sum_{U\in \cP}\left(\sum_{S\in \binom{V}{k-1}}\prob{A_{U,S}}+\sum_{v\in V}\prob{A_{U,v}}\right) < \sum_{U \in \cP}\left(\frac{1}{2L}+\frac{1}{2L}\right) = 1,
\end{align*}
so there exists a choice of $\cP$ such that $\cP$ is a $(\delta,\beta,L,1)$-partition for $G$.
\end{proof}

\subsection{Partition sequences and vortices}
When iteratively constructing partitions of vertex sets using Lemma~\ref{lem:vortexpartition}, we will work with sequences of successively refined partitions of the following form.
\begin{definition}
    Given a finite set $V$, say that a sequence $\cP_0,\ldots,\cP_t$ of partitions of $V$ is an \defn{$L$-partition sequence} if $k\mid |V|$, $\cP_0=\{V\}$, and for each $s \in [t]$ and $W \in \cP_{s-1}$, $\cP_s$ is a refinement of $\cP_{s-1}$ and $\cP_s[W]$ is a $k$-equitable $L$-partition of $W$.
\end{definition}

We characterize the sizes of sets in $L$-partition sequences using the following fact.
\begin{fact} \label{fact:floorceiling}
    Let $n,k,\ell \in \bN$.
    Let $a \in \{\lflr\frac{n}{k}\rflr,\lcl\frac{n}{k}\rcl\}$ and $b\in \{\lflr\frac{a}{\ell}\rflr,\lcl\frac{a}{\ell}\rcl\}$.
    Then $b\in \{\lflr\frac{n}{k\ell}\rflr,\lcl\frac{n}{k\ell}\rcl\}$.
\end{fact}
\begin{proof}
Since $n$ is an integer, let
\[ n = xk\ell + yk+z = x'k\ell - y'k - z', \]
where $x,y,z,x',y',z'$ are integers, $0\le y,y' < \ell$, and $0 \le z,z' < k$.
Then
\[ b \ge  \lflr \frac{\lflr \frac{n}{k}\rflr}{\ell} \rflr = \lflr \frac{\lflr \frac{xk\ell+yk+z}{k}\rflr}{\ell} \rflr = \lflr \frac{x\ell+y}{\ell}\rflr= x = \lflr \frac{n}{k\ell} \rflr. \]
Similarly,
\[ b \le \lcl\frac{\lcl\frac{x'k\ell-y'k -z'}{k}\rcl}{\ell}\rcl = \lcl\frac{x'\ell-y'}{\ell}\rcl = x' = \lcl \frac{n}{k\ell}\rcl.  \]
Thus, since $b\in \bN$, we have $b \in \left\{\lflr\frac{n}{k\ell}\rflr,\lcl\frac{n}{k\ell}\rcl\right\}$.
\end{proof}

The following lemma shows that all partitions in an $L$-partition sequence, induced on any of the partition sets, are $k$-equitable.
\begin{lemma} \label{lem:partitionsequence}
    Let $\cP_0,\ldots,\cP_t$ be an $L$-partition sequence for a finite set $V$.
    Then for all $0\le r \le s \le t$ and $W \in \cP_r$, $\cP_s[W]$ is a $k$-equitable $L^{s-r}$-partition of $W$.
\end{lemma}
\begin{proof}
Use induction on $q:=s-r$.
If $q=0$, then for $W \in \cP_r$ and $U \in \cP_s[W]$, we have $U=W$ and $k \mid |U|$, so we are done in this case.

Suppose $q \ge 1$, and let $W \in \cP_r$ and $U \in \cP_s[W]$.
Then there is a unique set $X \in \cP_{s-1}[W]$ such that $U \subseteq X \subseteq W$.
By the induction hypothesis, $|X| \in \left\{k\lflr\frac{|W|}{kL^{q-1}}\rflr,k\lcl\frac{|W|}{kL^{q-1}}\rcl\right\}$.
Because $\cP_s[X]$ is a $k$-equitable $L$-partition of $X$, $|U| \in \left\{k\lflr\frac{|X|}{kL}\rflr,k\lcl\frac{|X|}{kL}\rcl\right\}$.
Thus, by Fact~\ref{fact:floorceiling}, $|U|/k \in \left\{\ff{|W|}{kL^q},\cf{|W|}{kL^{q}}\right\}$, i.e., $|U| \in \left\{k\ff{|W|}{kL^q},k\cf{|W|}{kL^{q}}\right\}$.
Therefore, the lemma follows by induction.
\end{proof}
By applying Lemma~\ref{lem:partitionsequence} with $r=0$, we obtain the following.
\begin{cor} \label{cor:partitionsequence}
    Let $\cP_0,\ldots,\cP_t$ be an $L$-partition sequence for a finite set $V$, and let $n := |V|$ and $\alpha := 1/L$.
    Then for all $0\le s \le t$ and $U \in \cP_s$, we have $|U| \in \left\{k\lflr\frac{\alpha^s n}{k}\rflr,k\lcl\frac{\alpha^s n}{k}\rcl\right\}$.
\end{cor}

We now state the definition of a vortex, which is a sequence of refined $(\delta,\beta,L,1)$-partitions.
\begin{definition}
    Given a $k$-graph $G$, say that an $L$-partition sequence $\cP_0,\ldots,\cP_t$ of $V(G)$ with $t \ge 1$ is a \defn{$(\delta,\beta,L,C)$-vortex} for $G$ if, setting $\alpha := 1/L$, the following hold:
    \begin{enumerate}[label=(V\arabic*)]
        \item For all $1\le s \le t$ and $W\in \cP_{s-1}$, $\cP_s[W]$ is a $(\delta,\beta,L,1)$-partition for $G[W]$;
        \item $C\le \alpha^tn \le LC$.
    \end{enumerate}
\end{definition}

We are now ready to prove the main result of Section~\ref{sec:vortex}, showing the existence of vortices in vertex-regular Dirac graphs.
\begin{lemma}[Vortex Lemma] \label{lem:vortex}
    Let $k\ge 3$, let $1/C \ll 1/L \ll \delta$ and $n \ge LC$ with $k \mid n$.
    Suppose $G$ is an $n$-vertex $r$-regular $k$-graph with $\delta_{k-1}(G) \ge (\frac12+\delta)n$ for some $r$.
    Then there exists a $(\delta,r/n^{k-1},L,C)$-vortex $\cP_0,\ldots,\cP_t$ for $G$.
\end{lemma}
\begin{proof}[Proof of Lemma~\ref{lem:vortex}]
Define $\alpha := 1/L$, $t = t(n):= \lflr\frac{\log{n}-\log{C}}{\log{L}}\rflr$, and $\beta := r/n^{k-1}$, and note that $\beta = \Theta_k(1)$ and $C \le \alpha^t n \le LC$.
We will prove Lemma~\ref{lem:vortex} by induction on $t$.
If $t=1$, then by Lemma~\ref{lem:vortexpartition}, there exists a $(\delta,\beta,L,1)$-partition $\cP$ for $G$, so by setting $\cP_0 := \{V(G)\}$ and $\cP_1 :=\cP$, (V1) holds.

Suppose $t \ge 2$.
By the induction hypothesis, there exists a $(\delta,\beta,L,LC)$-vortex $\cP_0,\allowbreak \ldots,\allowbreak\cP_{t-1}$ for $G$.
Let $W \in \cP_{t-1}$.
By property (P3) of $W$, $G[W]$ is $|W|^{-1/3}$-almost $\beta |W|^{k-1}$-regular.
By property (P2), $\delta_{k-1}(G[W]) \ge (\frac12+\delta - |W|^{-1/3})|W|$.
Thus, by Lemma~\ref{lem:vortexpartition} applied to $G[W]$, there exists a $(\delta,\beta,L,1)$-partition $\cP_W$ for $G[W]$.

Let $\cP_t := \bigcup_{W \in \cP_{t-1}}\cP_W$, and consider the sequence $\cP_0,\ldots,\cP_t$ of partitions.
Then for all $W \in \cP_{t-1}$, $\cP_t[W] = \cP_W$ is a $(\delta,\beta,L,1)$-partition for $G[W]$, so (combined with the induction hypothesis) $\cP_0,\ldots,\cP_t$ is an $L$-partition sequence and (V1) holds for $s=t$.
(V2) follows from the definition of $t$, and (V1) for $s<t$ follow from the induction hypothesis.
Therefore, $\cP_0,\ldots,\cP_t$ is a $(\delta,\beta,L,C)$-vortex, and so the result follows by induction.
\end{proof}

\section{Cover-down lemma} \label{sec:coverdown}
We begin this section by proving Lemma~\ref{lem:coverdown}.
The overall proof strategy is as follows.

Our goal is to decompose all edges in $G-G[\cP]$, i.e., ``crossing" edges which are each not contained entirely in a set in the partition $\cP$, into perfect matchings by borrowing a small number of edges from $G[\cP] = \bigcup_{U\in \cP}G[U]$.
Because $G-G[\cP]$ is almost regular, one could apply Theorem~\ref{thm:approxdecomp} to decompose almost all edges of $G-G[\cP]$ into perfect matchings.
By defining an auxiliary graph $K$ whose vertices are the leftover edges in $G-G[\cP]$ not in the decomposition, drawing an edge between pairs with non-empty intersection, one could apply Theorem~\ref{thm:equitable} to find an equitable colouring of this graph, which corresponds to a decomposition of the leftover edges in $G-G[\cP]$ into approximately equal-sized matchings.
One would then hope to complete each of these matchings $M$ to a perfect matching by borrowing edges from $G[\cP]$.
However, if $|U\setminus V(M)|$ is not divisible by $k$ for some $U \in \cP$, then we cannot cover all of $U$ when adding edges in $G[U]$ to $M$ to create a larger matching.

We fix this divisibility issue by focusing only on edges in $G-G[\cP]$ with a fixed intersection size with each $U \in \cP$ and applying Corollary~\ref{cor:kequitable} to an auxiliary graph $K_J \subseteq K$ to decompose the edges into matchings $M$ for which $|U\setminus V(M)|$ is divisible by $k$ for all $U\in \cP$.
We can then apply Lemma~\ref{lem:sparsematchings} to find a matching on $U\setminus V(M)$ in $G$ for each $U \in \cP$, which we borrow from $G[\cP]$ to complete $M$ to a perfect matching.
However, this creates another problem: we have no control over the structure of $K_J$, so it may fail to satisfy the condition $|V(K_J)| > k \cdot \Delta(K_J)$ in Corollary~\ref{cor:kequitable}.
This can be remedied by using Lemma~\ref{lem:uniformtypes} to set aside a $k$-graph $F$ at the start of the proof that contains many edges of every possible intersection size with the sets in $\cP$.
We add the edges of $F$ back when constructing the auxiliary graph $K_J$, ensuring that the condition $|V(K_J)| > k \cdot \Delta(K_J)$ is satisfied for any choice of intersection sizes with the partition sets in $\cP$.
The proof then proceeds as expected.

\begin{proof}[Proof of Lemma~\ref{lem:coverdown}]
Define $\alpha :=1/L$ and choose new constants $\varrho,\mu,\eps_0,\eps_1,\eps_2,\delta_1,\delta_2$ such that
\[0 < \varrho \ll \mu \ll \eps_0 \ll \eps_1 \ll \eps_2 \ll \delta_1 \ll \delta_2 \ll \eta \ll \alpha \ll \delta <1.\]

Let $V:=V(G)$ and $\cP = \{U_1,\ldots,U_L\}$.
Let $\cI = \cI_{L,k}$ and $\cJ = \cJ_{L,k}$ as in Definition~\ref{def:types}, and note that $|\cI| \le L^k$.
Define $g:2^V \to \bN_0^L$ by $g(S) = (|S\cap U_1|,\ldots,|S\cap U_L|)$.
For each $J\in \cJ$ and an arbitrary subgraph $B \subseteq G$, we denote by $B_J$ the spanning subgraph of $B$ with $E(B_J) = E(B) \cap g^{-1}(J)$.

\begin{figure}
    \centering
    \includegraphics{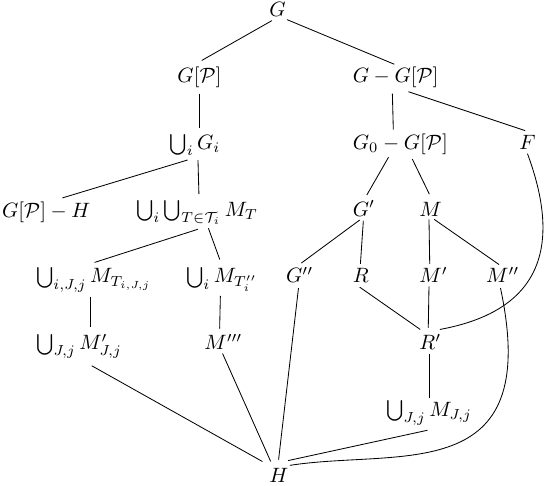}
    \caption{A diagram of the subgraphs of $G$ defined in Lemma~\ref{lem:coverdown}. Each edge of $G$ is in precisely the $k$-graphs along a unique maximal vertical path in the diagram.}
    \label{fig:coverdowninclusion}
\end{figure}

By properties (P1)--(P2) of the $(\delta,\beta,L,2)$-partition $\cP$ for $G$, $\cP$ satisfies the hypothesis of Lemma~\ref{lem:uniformtypes} with $\delta$ replaced by $\delta/2$, so there exists a $d$-regular spanning subgraph $F$ of $G-G[\cP]$ such that $\Delta_{k-1}(F) \le \eps_2 n$ and $|E(F_J)| \ge \eps_0 n^k$ for all $J\in \cJ$, where $d=\lflr \eps_1 n^{k-1}\rflr$.
Let $G_0 := G-F$, and observe that $d_{G_0,U}(S) \ge (\frac12+\delta-|U|^{-1/3})|U|-\eps_2n \ge (\frac12+\frac\delta2)|U|$, so $\cP$ is a $(\delta/2,\beta,L,2)$-partition for $G_0$.

\begin{claim}
    There exists a matching $M \subseteq G_0-G[\cP]$ such that for all $J \in \cJ$, $|E(M_J)|=k$.
\end{claim}
\begin{proof}[Proof of Claim]
We will show by induction that for all $0\le m \le k|\cJ|$, there exists a matching $M^{(m)} \subseteq G_0-G[\cP]$ such that $|E(M^{(m)})| = m$ and $|E(M^{(m)}_I)|\le k$ for all $I\in \cJ$.
If $m=0$, we are done by setting $M^{(0)}$ to be empty, so suppose $m \ge 1$, and let $M' := M^{(m-1)}$.
Since $m-1< k |\cJ|$, let $J\in \cJ$ such that $|E(M'_J)| < k$.
Let $j \in [L]$ such that $J_{j}>0$, and define $J' \in \bN_0^L$ by $J'_{j}=J_{j}-1$ and $J'_i = J_i$ for $i \neq j$.

For each $i \in [L]$, $|U_i \cap V(M')| < k^2|\cJ| \le k^2L^k \ll |U_i|$, so there exists $S_i \in \binom{U_i\setminus V(M')}{J'_i}$; let $S := \bigcup_{i=1}^LS_i \in \binom{V}{k-1}$.
Because $\cP$ is a $(\delta/2,\beta,L,2)$-partition for $G_0$, $d_{G_0,U_j}(S) \ge \frac12 |U_j|$, so there are at least $\frac12 |U_j| - k^2L^k>0$ choices for $v\in U_{j}\setminus V(M')$ such that $S\cup \{v\} \in E(G_0)$.

Choose such a vertex $v$, set $M^{(m)} := M' \cup \{S\cup \{v\}\}$, and note that $M^{(m)}$ is a matching in $G_0-G[\cP]$ with $|M^{(m)}_I| \le k$ for all $I\in \cJ$.
Therefore, by induction, we are done by setting $M:= M^{(k|\cJ|)}$.
\end{proof}
Let $r_0 := r-d$, $\beta_0 := r_0/n^{k-1}$, and $G' := G_0-G[\cP]-M$.
Note that $G_0$ is $\gamma$-almost $\beta_0n^{k-1}$-regular for $\gamma = \frac{2n^{-1/3}r}{r-d} =  O(n^{-1/3})$ (note that $r=\Omega_k(n^{k-1})$ and $d = O(\eps_1 n^{k-1})$ imply $r/(r-d)=O(1)$).
By the minimum $(k-1)$-degree conditions on $G_0$ and $G[U_1]$ (the latter via property (P2)), $\beta_0,\beta = \Theta_k(1)$.
More precisely, $(\frac12 + \frac\delta2) \binom{n}{k-1} \le \beta_0 n^{k-1} \le \binom{n-1}{k-1}$, so $\frac{1+\Omega(\delta)}{2(k-1)!} \le \beta_0 \le \frac{1}{(k-1)!}$, and similarly $\frac{1}{2(k-1)!} \le \beta \le \frac{1+o(1)}{(k-1)!}$, so $\beta_0 \le 2 \beta$ and $\beta \le 2 \beta_0$.

\begin{claim} \label{clm:iterativeapprox}
    Let $r' := (\beta_0-\alpha^{k-1}\beta)n^{k-1}$.
    Then $G'$ is $\varrho$-almost $r'$-regular and contains $\lflr (1-\mu)r'\rflr$ edge-disjoint perfect matchings.
\end{claim}
\begin{proof}[Proof of Claim]
Since $1/n \ll \mu \ll \alpha \ll \delta$ and $|U_i| = \alpha n \pm k$ for all $i$,
\begin{equation}
    \begin{aligned} \label{eq:iterativecodegree1}
        \delta_{k-1}(G') &\ge \delta_{k-1}(G) - \Delta_{k-1}(F)- \Delta_{k-1}(G[\cP])-\Delta_{k-1}(M) \\
        &\ge \left(\frac12+\delta\right)n- \eps_2 n - (\alpha n+k)-1 \\
        &\ge \left(\frac12+\mu\right)n.
    \end{aligned}
\end{equation}
By property (P3) of the $(\delta,\beta,L,2)$-partition $\cP$, for $v\in V$, since $|U_i| = \alpha n \pm k$ for all $i$,
\begin{align*}
    d_{G'}(v) &= d_{G_0}(v)-d_{G[\cP]}(v)-d_M(v) \\
    &= \left(1\pm O\left(n^{-1/3}\right)\right)\beta_0n^{k-1} - \left(1\pm 2(\alpha n - k)^{-1/3}\right)\beta (\alpha n \pm k)^{k-1}\pm 1 \\
    &= \left(1\pm O\left(n^{-1/3}\right)\right)(\beta_0-\alpha^{k-1} \beta)n^{k-1}.
\end{align*}
Thus, because $n^{-1/3} \ll \varrho \ll \mu$, $G'$ is $\varrho$-almost $r'$-regular, and the claim follows from Theorem~\ref{thm:approxdecomp}.
\end{proof}

Let $r_0' := (\beta - \alpha^{k-1}\beta)n^{k-1}$ and $r'' := \lflr r' - 3 \mu r_0' \rflr$.
Since $\beta_0 \le 2\beta$, we have $r' \le 3 r_0'$, and so $r'' \le (1-\mu)r'$.
Thus, let $G''$ be the spanning subgraph of $G'$ whose edge set is the union of $r''$ of the matchings from Claim~\ref{clm:iterativeapprox}, and let $R:= G' - G''$.
$R$ represents the remaining edges that have not been decomposed into perfect matchings, and our goal is to add back the edges in $F$ and borrow a small number of edges from $G[\cP]$ and $M$ to decompose $R$, plus these edges, into perfect matchings.
In order to guarantee that divisibility within each $U\in \cP$ is maintained when constructing perfect matchings using edges from $G[\cP]$, we will focus only on edges in $R_J \cup F_J \cup M_J$ for some $J\in \cJ$, which have fixed intersection sizes with each of sets in $\cP$.

For each $J\in \cJ$, let $0\le a_J <k$ with $a_J \equiv -|E(R_J \cup F_J)| \pmod{k}$.
Let $M'\cup M''$ be a partition of $M$ such that $|E(M_J')| = a_J$ for all $J\in \cJ$, and let $R' := R\cup F \cup M'$, so that $R'_J=R_J\cup F_J \cup M_J'$ for $J\in \cJ$.
Observe that
\begin{equation} \label{eq:iterativepartition}
    E(G) = E(G'') \sqcup E(R') \sqcup E(M'') \sqcup E(G[\cP]).
\end{equation}
For $J\in \cJ$, define the auxiliary graph $K_J$ by $V(K_J) = E(R'_J)$ and $E(K_J) = \{\{e,f\}:e\cap f \neq \emptyset\}$.
We will apply Corollary~\ref{cor:kequitable} to obtain a $k$-equitable colouring of $K_J$, which corresponds to a decomposition of $R_J'$ into approximately equal-sized matchings.
\begin{claim} \label{clm:iterativeequitable}
    $K_J$ has a $k$-equitable $\ell$-colouring for $\ell = \lflr \frac{\delta_1(\alpha n)^{k-1}}{2kL^k} \rflr$.
\end{claim}
\begin{proof}[Proof of Claim]
Because $\varrho \ll \mu$, and $\beta,\beta_0 = \Theta_k(1)$ implies $r_0' = O(r')$,
\[ \Delta_1(R'_J) \le \Delta_1(R) +\Delta_1(F)+\Delta_1(M'_J) \le O\left(\mu n^{k-1}\right)+\eps_1n^{k-1}+k = O\left(\eps_1n^{k-1}\right). \]
Because $\eps_1 \ll \delta_1 \ll \alpha$,
\[ \frac{\Delta_1(R_J')}{\ell} = O\left(\frac{kL^k\eps_1 n^{k-1}}{\delta_1(\alpha n)^{k-1}}\right) = O_k\left(\frac{\eps_1}{\delta_1 \alpha^{2k-1}}\right) < \frac{1}{k^2}. \]
For $e\in V(K_J)$, 
\[d_{K_J}(e) \le \sum_{v\in e}d_{R_J'}(v) \le k\cdot \Delta_1(R_J'), \]
and 
\[|E(R'_J)| = \frac{1}{k}\sum_{v\in V}d_{R_J'}(v) \le \frac{n}{k}\cdot \Delta_1(R_J'),\]
so it follows that
\[ k \cdot \Delta(K_J) \le k^2 \cdot \Delta_1(R_J') < \ell \]
and (as we will need later on)
\begin{equation} \label{eq:iterativefewrem}
    \frac{|E(R'_J)|}{\ell} \le \frac{n\cdot \Delta_1(R_J')}{k\ell} = O_k\left(\frac{\eps_1 n}{\delta_1\alpha^{2k-1}}\right) = O(\eps_2n).
\end{equation}
Further, since $F$ was obtained from Lemma~\ref{lem:uniformtypes},
\[ k\cdot \Delta(K_J) \le k^2 \Delta_1(R_J') = O_k\left(\eps_1 n^{k-1}\right) < \eps_0 n^k \le |E(F_J)| \le |E(R'_J)| = |V(K_J)|. \]
Thus, $\ell,|V(K_J)| > k \cdot \Delta(K_J)$ and $|V(K_J)| \equiv |E(R_J\cup F_J)|+a_J \equiv 0 \pmod{k}$, so by Corollary~\ref{cor:kequitable}, $K_J$ indeed has a $k$-equitable $\ell$-colouring.
\end{proof}
By Claim~\ref{clm:iterativeequitable}, let $\Phi_J = \{M_{J,1},\ldots,M_{J,\ell}\}$ be a $k$-equitable colouring of $K_J$.
By the definition of $K_J$, each $M_{J,j}$ is a matching in $R_J'$ of size $k\lflr\frac{|E(R'_J)|}{k\ell}\rflr$ or $k\lcl\frac{|E(R'_J)|}{k\ell}\rcl$.
For each $J\in \cJ$, $i \in [L]$, and $j\in [\ell]$, let $T_{i,J,j} := V(M_{J,j}) \cap U_i$.
For $i\in [L]$, let $T_i'' := V(M'')\cap U_i$ and
\[\cT_i := \{T_{i,J,j}:(J,j) \in \cJ \times [\ell]\}\cup \{T_i''\}.\]
We treat $\cT_i$ as a multiset consisting of $\ell|J|+1$ subsets of $U_i$.
Each element of $\cT_i$ consists of the vertices in $U_i$ covered by a small matching that we want to complete to a perfect matching using edges in $G[\cP]$.
The following claim details how these edges in $G[\cP]$ are chosen.
\begin{claim} \label{clm:iterativesparse}
    For each $i\in [L]$, there exist pairwise edge-disjoint matchings $(M_{T})_{T\in \cT_i}$ in $G[U_i]$ such that $V(M_{T}) = U_i\setminus T$ for all $T\in \cT_i$ and $\Delta_{k-1}(\bigcup_{T\in \cT_i}M_{T}) \le \delta_2|U_i|$.
\end{claim}
\begin{proof}[Proof of Claim]
Let $i\in [L]$.
First, observe that
\[|\cT_i| = \ell|\cJ|+1 \le \frac{\delta_1|U_i|^{k-1}}{kL^k}kL^k = \delta_1|U_i|^{k-1}.\]
Also, 
\[\delta_{k-1}(G[U_i]) \ge \left(\frac12+\delta - |U_i|^{-1/3}\right)|U_i| \ge \left(\frac12+\frac\delta2\right)|U_i|.\]
$(T_{i,J,j})_{(J,j) \in \cJ\times [\ell]}$ are subsets of $U_i$ such that, by \eqref{eq:iterativefewrem} and the fact that $\Phi_J$ is $k$-equitable,
\[|T_{i,J,j}| \le k|E(M_{J,j})| \le k\left(\frac{|E(R_J')|}{\ell}+k\right) = O_k(\eps_2 n)= O_k\left(\eps_2L|U_i|\right) \le \delta_1 |U_i|, \]
and $|T_{i,J,j}| \equiv J_i \cdot |E(M_{J,j})| \equiv 0 \pmod{k}$.
$T_i''$ is a subset of $U_i$ such that $|T_i''| \le k^2L^k \le \delta_1|U_i|$, and by \eqref{eq:iterativepartition}, since $k\mid |U_i|$, $k \mid |E(R_J')|$ for $J \in \cJ$, and $k\mid \sum_{v\in U_i}d_G(v)$ by the hypothesis of the lemma,
\begin{align*}
    |T_i''| &= \sum_{v\in U_i}d_{M''}(v) = \sum_{v\in U_i}(d_{G}(v)-d_{G''}(v)-d_{R'}(v)-d_{G[\cP]}(v)) \\
    &= \sum_{v\in U_i}d_G(v)- r''|U_i|-\sum_{J\in \cJ}J_i|E(R'_J)| - k|E(G[U_i])| \equiv 0 \pmod{k}.
\end{align*}
    
Therefore, the claim follows from Lemma~\ref{lem:sparsematchings} with $\delta$ replaced by $\delta/2$.
\end{proof}
For each $J\in \cJ$ and $j \in [\ell]$, let $M_{J,j}' := \bigcup_{i=1}^LM_{T_{i,J,j}}$, and let $M''' := \bigcup_{i=1}^L M_{T_i''}$, where the $M_{T_{i,J,j}}$ and $M_{T_i''}$ are from Claim~\ref{clm:iterativesparse}.
Observe that $M_{J,j} \cup M_{J,j}'$ is a perfect matching on $V$ for all $J\in \cJ$ and $j \in [\ell]$, and $M'' \cup M'''$ is a perfect matching on $V$ (since by Claim~\ref{clm:iterativesparse}, each of these matchings covers every $U_i$).
Let 
\[
    H := G'' \cup \left( \bigcup_{J\in \cJ}\bigcup_{j=1}^\ell (M_{J,j}\cup M_{J,j}')\right)\cup (M''\cup M''').
\]
We claim that $H$ satisfies the conclusion of the lemma.
Indeed, by \eqref{eq:iterativepartition},  $G-G[\cP] \subseteq H$ (note that $R' = \bigcup_{J,j}M_{J,j}$).
Further, since $H[\cP] = \bigcup_{i=1}^L\bigcup_{T\in \cT_i} M_{T}$ and $\delta_2 \ll \eta$, by Claim~\ref{clm:iterativesparse}, we have
\[ \Delta_{k-1}(H[\cP]) \le \delta_2 \max_i|U_i| \le \eta n. \]
Finally, note that $G''$ has a decomposition into $r''$ perfect matchings, and $M_{J,j} \cup M_{J,j}'$ ($J \in \cJ$, $j \in [\ell]$) and $M'' \cup M'''$ are each perfect matchings, so $H$ admits a decomposition into perfect matchings.
\end{proof}

We next observe that a decomposition of a $k$-graph $G$ into matchings which are perfect on vertex sets in a sequence of successively refined partitions can be easily transformed into a $1$-factorization of $G$.
For a vertex set $V$ and $\cS \subseteq 2^V$, an \defn{$\cS$-matching} is a matching $M$ with $V(M) \in \cS$.
\begin{lemma} \label{lem:S-matching}
    Let $\cP_0,\ldots,\cP_t$ be partitions of a finite set $V$ such that $\cP_0 = \{V\}$ and $\cP_s$ is a refinement of $\cP_{s-1}$ for all $s\in [t]$, and let $\cS := \bigcup_{s=0}^t\cP_s$.
    Suppose that an $r$-regular $k$-graph $G$ on the vertex set $V$ admits a decomposition into $\cS$-matchings.
    Then $G$ admits a decomposition into $r$ perfect matchings.
\end{lemma}
\begin{proof}
Let $\cM := \{M_1,\ldots,M_d\}$ be the decomposition of $G$ into $\cS$-matchings.
We will use induction on $t+d$.
If $t=0$, then since $\cP_0 = \{V\}$, $\cM$ itself is a $1$-factorization of $G$, and if $d=0$, then $G$ is empty and the result is trivial.

Suppose $t,d\ge 1$.
Let $r$ be the vertex degree of $G$.
If there exists a matching $M \in \cM$ with $V(M) = V$, then $G-M$ has a decomposition $\cM\setminus \{M\}$ into $\cS$-matchings and is $(r-1)$-regular, and so $G-M$ has a $1$-factorization $M_1',\ldots,M_{r-1}'$ by the induction hypothesis.
Thus, $M,M_1',\ldots,M_{r-1}'$ is a decomposition of $G$ into $r$ perfect matchings.

Suppose otherwise, and note that $V(M) \in \bigcup_{s=1}^t\cP_s$ for all $M\in \cM$, so $G=G[\cP_1]$.
Thus, for $v\in U \in \cP_1$, $d_{G[U]}(v) = d_G(v) = r$, so $G[U]$ is $r$-regular.
Now, by the induction hypothesis applied to the $k$-graph $G[U]$, the sequence of partitions $\cP_0',\ldots,\cP_{t-1}'$ with $\cP_s' = \cP_{s+1}[U]$ for all $s$, and the matchings $\{M\in \cM:V(M)\subseteq U\}$, $G[U]$ admits a $1$-factorization $M_{U,1},\ldots,M_{U,r}$.
Thus, $G$ admits a $1$-factorization $(\bigcup_{U\in \cP_1}M_{U,1}),\ldots,(\bigcup_{U\in \cP_1}M_{U,r})$.
Therefore, the lemma follows by induction.
\end{proof}

We will now apply Lemma~\ref{lem:coverdown} iteratively to obtain the following strong approximate decomposition result, which decomposes almost all edges of $G$ into matchings which are perfect on vertex sets in vortex partitions, leaving only edges in $G[\cP_t]$, where $\cP_t$ is a partition of $V(G)$ into constant-size sets.
\begin{lemma} \label{lem:iterative}
    Let $k \ge 3$, $1/n \le 1/C \ll 1/L \ll \delta < 1$ with $k \mid n$, $t \ge 0$, and $\alpha := 1/L$.
    Suppose $G$ is an $r$-regular $n$-vertex $k$-graph for some $r$.
    Let $\cP_0,\ldots,\cP_{t+1}$ be a $(\delta,\beta,L,C)$-vortex for $G$ for some $\beta$, and let $\cS := \bigcup_{s=0}^t \cP_s$.
    Then there exists a $k$-graph $H$ such that $G-G[\cP_t] \subseteq H \subseteq G-G[\cP_{t+1}]$, $\Delta_{k-1}(H[\cP_t]) \le \alpha^{t+2}n$, and $H$ admits a decomposition into $\cS$-matchings.
\end{lemma}
\begin{proof}
Let $\eta$ be a new constant with $1/C \ll \eta \ll 1/L$.
Use induction on $t$.
If $t=0$, we are done by setting $H$ to be empty.

Suppose $t\ge 1$.
By the induction hypothesis, since $\cP_0,\ldots,\cP_t$ is a $(\delta,\beta,L,LC)$-vortex for $G$, let $H'$ be a $k$-graph such that $G-G[\cP_{t-1}] \subseteq H' \subseteq G-G[\cP_t]$, $\Delta_{k-1}(H'[\cP_{t-1}]) \le \alpha^{t+1}n$, and $H'$ admits a decomposition $\cM'$ into $\cS'$-matchings for $\cS' = \bigcup_{s=0}^{t-1}\cP_s$.
Let $W \in \cP_{t-1}$ and $G_W := G[W]-H'- G[\cP_{t+1}]$.
We would like to show that $G_W$ and $\cP_t[W]$ satisfy the hypotheses of Lemma~\ref{lem:coverdown}.

Letting $r' := |\{M \in \cM':W \subseteq V(M)\}|$, we have that $d_{G-H'}(v) = r-r'$ for all $v\in W$.
Since $G-H' \subseteq G[\cP_{t-1}]$, this means that $G[W]-H'$ is $(r-r')$-regular; let $\beta' := (r-r')/|W|^{k-1}$.
Also, for $U\in \cP_t[W]$ and $S \in \binom{W}{k-1}$, since $\Delta_{k-1}(H'[W]) \le \alpha^{t+1}n \le 2\alpha |U|$ and $\Delta_{k-1}(G[\cP_{t+1}]) \le \max_{P \in \cP_{t+1}}|P| \le \alpha |U|+k$,
\[ d_{G_W,U}(S) \ge \left(\frac12+\delta-|U|^{-1/3}\right)|U|-\alpha^{t+1}n-\alpha |U|-k \ge \left(\frac12+\frac\delta2\right)|U|. \]
By summing over all $U$, we obtain $\delta_{k-1}(G_W) \ge \left(\frac12+\frac\delta2\right)|W|$, so in particular, $\beta' = \Theta(1)$.
For $U\in \cP_t[W]$ and $v\in P \in \cP_{t+1}[U]$, since $|P| = \alpha^2 |W| \pm k$ and $|P| = \alpha |U|\pm k$,
\[ d_{G_W}(v) = \beta'|W|^{k-1} - \left(1\pm |P|^{-1/3}\right)\beta |P|^{k-1} = \left(1\pm |W|^{-1/3}\right)\left(\beta'-\alpha^{2(k-1)}\beta\right)|W|^{k-1} \]
and
\[ d_{G_W[U]}(v) = \left(1\pm |U|^{-1/3}\right)\beta|U|^{k-1}-\left(1 \pm |P|^{-1/3}\right)\beta|P|^{k-1} = \left(1\pm 2|U|^{-1/3}\right)\left(1-\alpha^{k-1}\right)\beta|U|^{k-1}, \]
so $G_W$ is $|W|^{-1/3}$-almost vertex-regular and $G_W[U]$ is $2|U|^{-1/3}$-almost $(1-\alpha^{k-1})\beta|U|^{k-1}$-regular.
Thus, $\cP_t[W]$ is a $(\delta/2,(1-\alpha^{k-1})\beta,L,2)$-partition for $G_W$.
Finally, for $U \in \cP_t[W]$, 
\[\sum_{v\in U}d_{G_W}(v) = (r-r')|U|-k|E(G[\cP_{t+1}[U]])|  \equiv 0 \pmod{k}. \]
Therefore, by Lemma~\ref{lem:coverdown} with $\delta$ replaced by $\delta/2$, there exists a subgraph $H_W$ of $G_W$ such that $G_W-G_W[\cP_t] \subseteq H_W$, $\Delta_{k-1}(H_W[\cP_t]) \le \eta |W|$, and $H_W$ admits a $1$-factorization.

Let $H := H' \cup \bigcup_{W\in \cP_{t-1}}H_W$.
Since $G-G[\cP_{t-1}] \subseteq H' \subseteq G-G[\cP_t]$ and $G[W]- H' -G[\cP_t] \subseteq H_W \subseteq G-G[\cP_{t+1}]$ for all $W\in \cP_{t-1}$, we have $G-G[\cP_t] \subseteq H \subseteq G-G[\cP_{t+1}]$.
Also, for $W \in \cP_{t-1}$, $U \in \cP_t[W]$, and $S\in \binom{U}{k-1}$, we have 
\[d_{H[U]}(S) = d_{H_W[U]}(S) \le \eta |W| \le \eta (\alpha^{t-1}n+k) \le \alpha^{t+2}n,\]
so $\Delta_{k-1}(H[\cP_t]) \le \alpha^{t+2}n$.
Finally, because each $H_W$ was defined within $G_W$ and is thus edge-disjoint from $H'$, and perfect matchings in $H_W$ are $\cS$-matchings in $H$, $H$ admits a decomposition into $\cS$-matchings.
This completes the proof of Lemma~\ref{lem:iterative}.
\end{proof}

\section{Final absorption lemma} \label{sec:absorption}
The goal of this section is to prove Lemma~\ref{lem:absorbing}, which shows the existence of an absorbing structure $A$ that can be set aside before performing the iterative absorption process in Section~\ref{sec:coverdown} that absorbs any resulting leftover.
Recall from Section~\ref{sec:outline} that we construct the absorbing structure $A$ through a series of reductions, which we present in the opposite order as a multi-step process in which an arbitrary leftover is absorbed by fixed absorbers into structures with increasingly strong properties, eventually yielding a $1$-factorization of the leftover along with $A$.

\subsection{Random vertex sets}
The following lemma shows that large sets chosen (almost) uniformly at random from the vertices of a $k$-graph $G$ are likely to approximately preserve a minimum $(k-1)$-degree condition of $G$. This is a well-known fact; we include the short proof for completeness.

\begin{lemma} \label{lem:randomset}
    Let $k\ge 3$, let $n,q,z \in \bN$ and $\eps \in (0,1)$ such that $z+q\le n$ and $z/n \le \eps/4$, let $G$ be an $n$-vertex $k$-graph with $\delta_{k-1}(G) \ge (\frac12 + \eps)n$, and let $Z \subseteq V(G)$ with $|Z| \le z$.
    Let $Q$ be a uniformly random $q$-subset of $V(G)\setminus Z$, and let $S \in \binom{V}{k-1}$.
    Then
    \[\prob{d_{G,Q}(S) \le \left(\frac12+\frac\eps2\right)q} \le e^{-\Omega(\eps^2q)}. \]
\end{lemma}

\begin{proof}
Let $V:= V(G)$ and $Y_S := d_{G,Q}(S)$.
Observe that
\[Y_S \sim \op{Hypergeometric}(n-|Z|,d_{G,V\setminus Z}(S),q),\]
since $Y_S$ is obtained by choosing a uniformly random $q$-subset $Q$ of the $(n-|Z|)$-set $V\setminus Z$ and counting the number of ``successes" $v\in V\setminus Z$ with $S\cup \{v\} \in E(G)$ that are contained in $Q$.
Since 
\[d_{G,V\setminus Z}(S) = d_G(S)-d_{G,Z}(S) \ge \left(\frac12+\eps\right)n-z \]
and $z/n \le \eps/4$,  we have that
\[ \expn{Y_S} = \frac{q\cdot d_{G,V\setminus Z}(S)}{n-|Z|} \ge \frac{q(\frac12+\eps-\frac{z}{n})n}{n} \ge \left(\frac12+\frac{3\eps}{4}\right)q \]
and $\expn{Y_S} \le q$.
Thus, by Theorem~\ref{thm:chernoffhyper},
\[ \prob{Y_S \le \left(\frac12+\frac\eps2\right)q} \le \prob{Y_S \le \left(1-\frac\eps4\right)\expn{Y_S}} \le e^{-\Omega(\eps^2q)}, \]
as desired.
\end{proof}

\subsection{3-step absorption process with compatible matchings}
Given a set $U_0$ in the final vortex partition, we will construct an absorber in three steps for $U_0$ which, when combined with any leftover in $U_0$, admits a $1$-factorization.
The first of these steps constructs an absorber $A_1$ that, together with any leftover in $U_0$, decomposes into matchings which each induce a union of sets in $\cP$ for some partition $\cP$ of $V(A_1)$.
This absorber is a multihypergraph containing ``fake" edges that do not actually exist in our original $k$-graph.

For multisets (resp. $k$-multigraphs) $S,T$, let $S+T$ denote the multiset (resp. $k$-multigraph) with multiplicities (resp. edge multiplicities) added between $S$ and $T$.
\begin{definition}
    Given a (multi)hypergraph $G$ and a collection $\cP$ of vertex sets, say that $G$ is \defn{$\cP$-compatible} if for $u,v \in V(G)$ and $U \in \cP$, $u,v \in U$ implies $d_G(u) = d_G(v)$.
    Say that a vertex set $V$ is \defn{$\cP$-compatible} if $U \cap V \in \{\emptyset, U\}$ for all $U \in \cP$.
    Say that $G$ is \defn{$\cP$-decomposable} if $G$ admits a decomposition into matchings, each of which is $\cP$-compatible.
\end{definition}
Note that if $\cP[V(G)]$ is a partition of $V(G)$, then a matching $M$ in $G$, treated as a subgraph of $G$, is $\cP$-compatible if and only if it induces a $\cP$-compatible vertex set $V(M)$.
Also, if $G$ is $\cP$-decomposable, then it is $\cP$-compatible.

The following lemma details the first step of the construction of the absorber.
\begin{lemma} \label{lem:absorbingmulti}
    Let $k\ge 3$.
    Let $m_0,m_1 \in k\bN$ such that $m_1 \ge m_0+\binom{m_0}{k}k(k-1)$.
    Let $U_0 \subseteq U_1$ be sets with $|U_0|=m_0$ and $|U_1|=m_1$, and let $\cP$ be a partition of $U_1$ into $k$-sets such that $U_0$ is $\cP$-compatible.
    Then there exists a $k$-multigraph $A_1$ with $V(A_1) = U_1$ and $|E(A_1)| = O\left(m_0^{k^2}\right)$ such that for any $\cP$-compatible $k$-graph $H_0$ with $V(H_0) = U_0$, $A_1 + H_0$ is $\cP$-decomposable.
\end{lemma}
\begin{proof}
Let $\overline{U_0} := U_1\setminus U_0$, and let $\cQ = \{Q_e:e \in \binom{U_0}{k}\}$ be a collection of pairwise disjoint subsets of $\overline{U_0}$ such that each $Q_e$ is the union of $k-1$ distinct sets in $\cP$, i.e., $\cQ$ is a partial partition of $\overline{U_0}$ into $\binom{m_0}{k}$ many $\cP$-compatible $k(k-1)$-sets.
This exists since $m_1 \ge m_0 + \binom{m_0}{k}k(k-1)$.

For a collection $\cS$ of subsets of $U_1$, denote $V(\cS) := \bigcup_{S\in \cS}S$.

\begin{figure}
    \centering
    \includegraphics[scale=0.8]{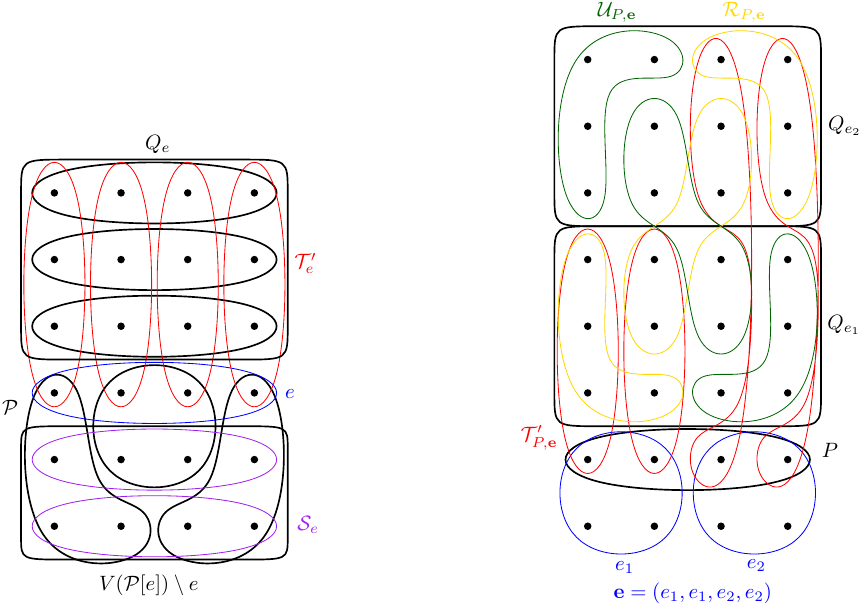}
    \caption{Examples of the matchings $\cS_e$, $\cT_e'$, $\cT'_{P,\mathbf{e}}$, $\cR_{P, \mathbf{e}}$, $\cU_{P,\mathbf{e}}$ that make up $A_1$ in the case $k=4$. If the edges $e,e_1,e_2$ depicted are present in $H_0$ (and $\mathbf{e} \in \cE_P'$), we form $\cP$-compatible matchings using $\cS_e$, $\cT'_{P,\mathbf{e}}$, and $\cU_{P,\mathbf{e}}$ and put $\cR_{P,\mathbf{e}}$ in $\cR$; if they are not present, we use the matchings $\cS_e$, $\cT_e'$, $\cU_{P,\mathbf{e}}$, and $\cR_{P,\mathbf{e}}$.}
    \label{fig:absorptionmulti}
\end{figure}

Let $e = \{v_1,\ldots,v_k\}\in \binom{U_0}{k}$, and note that $\cP[e] = \{P \in \cP: P\cap e \neq \emptyset\}$ is a collection of $x_e$ many $k$-sets for some $1\le x_e \le k$.
Let $\cS_e = \{S_{e,1},\ldots,S_{e,x_e-1}\}$ be a partition of $V(\cP[e])\setminus e$ into $k$-sets (since $|V(\cP[e])|-|e|=k(x_e-1)$), and let $\cT_e = \{T_{e,v_1},\ldots,T_{e,v_k}\}$ be a partition of $Q_e$ into $(k-1)$-sets (since $|Q_e| = k(k-1)$).
Let $\cT_e' = \{T_{e,v_1}',\ldots,T_{e,v_k}'\}$, where $T_{e,v}' := T_{e,v}\cup \{v\}$ for all $v\in e$.
\begin{claim} \label{clm:absorptionmulti1}
    For all $e \in \binom{U_0}{k}$, $\{e\}\cup \cS_e$ and $\cS_e \cup \cT_e'$ are $\cP$-compatible matchings.
\end{claim}
\begin{proof}[Proof of Claim]
We have $V(\{e\}\cup \cS_e) = V(\cP[e])$ and $V(\cS_e \cup \cT_e') = V(\cP[e])\cup Q_e$ by the definitions of $\cS_e$ and $\cT_e'$, and both of these vertex sets are $\cP$-compatible.
Also, $\cS_e$ is vertex-disjoint from $e$ and $\cT'_{e,v}$, since $V(\cS_e) = V(\cP[e])\setminus e \subseteq U_0\setminus e$ and $V(\cT'_{e,v}) = Q_e \cup e\subseteq \overline{U_0} \cup e$.
\end{proof}

For $v\in U_0$, let $E_v = \{e\in \binom{U_0}{k}:v\in e\}$.
Let $P \in \cP[U_0]$, $\cE_P := \prod_{v\in P}E_v$ (which consists of $k$-tuples of edges indexed by vertices in $P$), and $\mathbf{e} = (e_v)_{v\in P} \in \cE_P$.
Observe that $T_{e_u,u} \cap T_{e_{v},v} = \emptyset$ for distinct $u,v\in P$ ($\cQ$ and $\cT_e$ are partitions), so $\cT_{P,\mathbf{e}} := \{T_{e_v,v}:v\in P\}$ consists of pairwise disjoint $(k-1)$-subsets of $\overline{U_0}$, and $\cT_{P,\mathbf{e}}'=\{T_{e_v,v}':v\in P\}$ is a matching that covers $P$.
Let $V_{P,\mathbf{e}} := V(\cT_{P,\mathbf{e}})$.

Let $\cR_{P,\mathbf{e}}$ be a partition of $V_{P,\mathbf{e}}$ into $k$-sets (since $|V_{P,\mathbf{e}}|=k(k-1)$ is divisible by $k$).
Note that $\cQ[V_{P,\mathbf{e}}]$ is a collection of $y_\mathbf{e}$ many $k(k-1)$-sets for some $1\le y_\mathbf{e} \le k$.
Let $\cU_{P,\mathbf{e}}$ be a partition of $V(\cQ[V_{P,\mathbf{e}}])\setminus V_{P,\mathbf{e}}$ into $k$-sets (since $k$ divides $|V(\cQ[V_{P,\mathbf{e}}])| - |V_{P,\mathbf{e}}| = k(k-1)(y_\mathbf{e}-1)$).
\begin{claim} \label{clm:absorptionmulti2}
    For all $P \in \cP[U_0]$ and $\mathbf{e} \in \cE_P$, $\cT_{P,\mathbf{e}}'\cup \cU_{P,\mathbf{e}}$ and $\cR_{P,\mathbf{e}} \cup \cU_{P,\mathbf{e}}$ are $\cP$-compatible matchings.
\end{claim}
\begin{proof}[Proof of Claim]
We have $V(\cT_{P,\mathbf{e}}'\cup \cU_{P,\mathbf{e}}) = P \cup V(\cQ[V_{P,\mathbf{e}}])$ and $V(\cR_{P,\mathbf{e}} \cup \cU_{P,\mathbf{e}}) = V(\cQ[V_{P,\mathbf{e}}])$ by the definitions of $\cT_{P,\mathbf{e}}'$, $\cU_{P,\mathbf{e}}$, and $\cR_{P,\mathbf{e}}$, and both of these vertex sets are $\cP$-compatible.
Also, $\cU_{P,\mathbf{e}}$ is vertex-disjoint from $\cT_{P,\mathbf{e}}$ and $\cR_{P,\mathbf{e}}$, since $V(\cU_{P,\mathbf{e}}) = V(\cQ[V_{P,\mathbf{e}}])\setminus V_{P,\mathbf{e}} \subseteq \overline{U_0} \setminus V_{P,\mathbf{e}}$, $V(\cT_{P,\mathbf{e}}') = V_{P,\mathbf{e}}\cup P \subseteq V_{P,\mathbf{e}} \cup U_0$, and $V(\cR_{P,\mathbf{e}}) = V_{P,\mathbf{e}}$.
\end{proof}

Let $A_1$ be the $k$-multigraph with $V(A_1) = U_1$ and
\[ E(A_1) = \sum_{e\in \binom{U_0}{k}}(\cS_e \cup \cT_e') + \sum_{P \in \cP[U_0]}\sum_{\mathbf{e} \in \cE_P}(\cR_{P,\mathbf{e}}\cup \cU_{P,\mathbf{e}}). \]
Note that
\[ |E(A_1)| \le \binom{m_0}{k}((k-1)+k)+\frac{m_0}{k}\binom{m_0-1}{k-1}^k((k-1)+(k-1)^2) = O\left(m_0^{k^2}\right). \]

Let $H_0$ be a $\cP$-compatible $k$-graph with $V(H_0) = U_0$.
Let $\overline{H_0}$ be the complement $k$-graph of $H_0$, i.e., $V(\overline{H_0}) = U_0$ and $E(\overline{H_0}) = \binom{U_0}{k}\setminus E(H_0)$.

For $v\in P \in \cP[U_0]$, let $d_P := d_{H_0}(v)$ and $d_P' := d_{\overline{H_0}}(v)$ (which depend only on $P$ since $H_0$ is $\cP$-compatible), let $E_v \cap E(H_0) = \{e_{v,1},\ldots,e_{v,d_P}\}$, and let $E_v \cap E(\overline{H_0}) = \{e_{v,1}',\ldots,e_{v,d_P'}'\}$.
For $i\in [d_P]$, let $\mathbf{e}_{P,i} = (e_{v,i})_{v\in P} \in \cE_P$, and let $\cE_P' := \{\mathbf{e}_{P,i}: i \in [d_P]\}$.

\begin{claim} \label{clm:absorptionmulti3}
    $\cR := \bigcup_{P \in \cP[U_0]}\bigcup_{\mathbf{e}\in \cE_P'}\cR_{P,\mathbf{e}}$ is a $\cP$-compatible matching.
\end{claim}
\begin{proof}
Consider the vertex multiset
\begin{align*}
    V &:= \sum_{P \in \cP[U_0]}\sum_{\mathbf{e}\in \cE_P'}V_{P,\mathbf{e}} = \sum_{P \in \cP[U_0]}\sum_{v\in P}\sum_{i=1}^{d_P}V(T_{e_{v,i},v}) = \sum_{P\in \cP[U_0]}\sum_{v\in P} \sum_{e\in E_v \cap E(H_0)}V(T_{e,v}) \\
    &= \sum_{e \in E(H_0)}\sum_{v\in e}V(T_{e,v}) = \sum_{e\in E(H_0)} Q_e.
\end{align*}
Note that $\sum_{e\in E(H_0)}Q_e$ is a sum of pairwise disjoint $\cP$-compatible sets and is hence itself a $\cP$-compatible set.
Thus, $V$ is a $\cP$-compatible set (with multiplicities at most $1$).
Since each $\cR_{P,\mathbf{e}}$ is a matching on the vertex set $V_{P,\mathbf{e}}$, it follows that $\cR$ is a $\cP$-compatible matching.
\end{proof}

We can now decompose the edge multiset of $A_1+H_0$ as
\begin{align*}
    E(A_1 + H_0) &= \sum_{e\in E(H_0)}(\{e\}\cup \cS_e) + \sum_{e\in E(\overline{H_0})}(\cS_e \cup \cT_e') + \cR \\
    &+ \sum_{P \in \cP[U_0]}\left(\sum_{\mathbf{e} \in \cE_P'}(\cT_{P,\mathbf{e}}'\cup \cU_{P,\mathbf{e}})+ \sum_{\mathbf{e} \in \cE_P\setminus \cE_P'}(\cR_{P,\mathbf{e}}\cup \cU_{P,\mathbf{e}})\right).
\end{align*}
By Claims \ref{clm:absorptionmulti1}, \ref{clm:absorptionmulti2}, and \ref{clm:absorptionmulti3}, for all $e\in \binom{U_0}{k}$, $P\in \cP$, and $\mathbf{e} \in \cE_P$, the matchings $\{e\}\cup \cS_e$, $\cS_e \cup \cT_e'$, $\cR$, $\cT_{P,\mathbf{e}}'\cup \cU_{P,\mathbf{e}}$, and $\cR_{P,\mathbf{e}}\cup \cU_{P,\mathbf{e}}$ are $\cP$-compatible, so this gives a decomposition of $A_1+H_0$ into $\cP$-compatible matchings, as desired.
\end{proof}

The second step in the construction of our absorbing structure converts each ``fake" edge in the previous step into a matching that functions identically to the corresponding ``fake" edge, thus absorbing any leftover into matchings with the same properties as in Lemma~\ref{lem:absorbingmulti} but that actually exist in our original $k$-graph $G$.

\begin{lemma} \label{lem:absorbingcompatible}
    Let $k \ge 3$.
    Let $\eps > 0$ and $m_1,m_2 \in k \bN$ with $1/m_2 \ll \eps <1$ and $m_1/m_2 \le \eps/4$.
    Let $U_1 \subseteq U_2$ be sets with $|U_1|=m_1$ and $|U_2|=m_2$, and let $\cP_1$ be a partition of $U_1$ into $k$-sets.
    Let $G_2$ be a $k$-graph with $V(G_2) = U_2$ and $\delta_{k-1}(G_2) \ge (\frac12 + \eps)m_2$, and let $A_1$ be a $k$-multigraph with $V(A_1)=U_1$ and $|E(A_1)| \le \sqrt{m_2}$.
    Then there exists a partition $\cP_2$ of $U_2\setminus U_1$ into $k$-sets and a spanning subgraph $A_2$ of $G_2$ such that $U_1$ is an independent set in $A_2$, and for any $\cP_1$-compatible $k$-graph $H_0$ with $V(H_0)\subseteq U_1$, if $A_1 + H_0$ is $\cP_1$-decomposable, then $A_2 \cup H_0$ is $(\cP_1 \cup \cP_2)$-decomposable.
\end{lemma}
\begin{proof}
If $A_1$ is empty, then by setting $A_2$ to also be empty, the conclusion holds automatically, so suppose $A_1$ is non-empty.
Let $X$ be a $|E(A_1)|$-set disjoint from $\binom{U_1}{k}$, and let $f:X\to \binom{U_1}{k}$ be such that $|f^{-1}(e)|$ is the multiplicity of $e$ in $A_1$ for all $e\in E(A_1)$.
One can consider $X$ as a set of ``fake" edges and $f$ as the function that maps each fake edge to the $k$-set that it represents.

\begin{figure}
    \centering
    \includegraphics{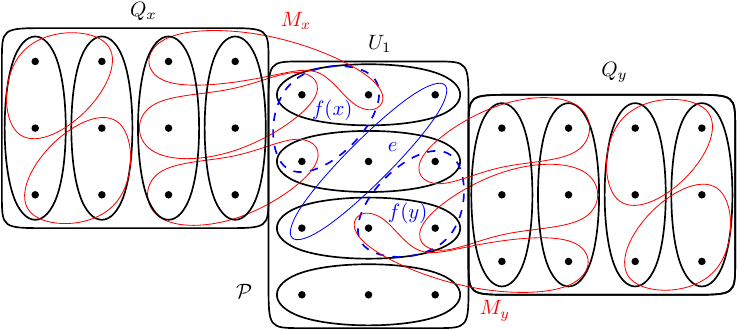}
    \caption{An example of a matching $M' = \{e\}\cup M_{x}\cup M_{y}$ (solid red and blue) in the decomposition of $H_0\cup A_2$ into $\cP$-compatible matchings in Lemma~\ref{lem:absorbingcompatible} in the case $k=3$.
    $M'$ is derived from the matching $M$ that consists of the edge $e \in E(H_0)$ and the ``fake edges" $x,y \in X$ corresponding to the edges $f(x),f(y) \in E(A_1)$ (dashed blue).}
    \label{fig:absorbingcompatible}
\end{figure}

Let $\cQ = \{Q_x:x\in X\}$ be a random $k$-equitable $|X|$-partition of $U_2\setminus U_1$ chosen such that each $Q_x$ is a uniformly random subset of $U_2\setminus U_1$ of a fixed size.
Since $m_1/m_2 \le \eps/4$ and $|Q_x| = \frac{m_2-m_1}{|X|}\pm k= \Omega(\sqrt{m_2})$ for all $x\in X$, by Lemma~\ref{lem:randomset} with $Z=U_1$ and a union bound, the probability that $d_{G_2,Q_x}(S) \le (\frac12+\frac\eps2)|Q_x|$ for some $x\in X$ and $S\in \binom{U_2}{k-1}$ is at most 
\[|X| \binom{m_2}{k-1}e^{-\Omega\left(\eps^2\sqrt{m_2}\right)} \le m_2^{k-1/2}e^{-\Omega\left(\eps^2m_2^{1/2}\right)} < 1, \]
so fix a choice of $\cQ$ such that, letting $G_x$ be $G_2[Q_x \cup f(x)]$ with the edge $f(x)$ removed (if it exists), we have $\delta_{k-1}(G_x) \ge (\frac12 +\frac\eps2)|Q_x| \ge (\frac12 + \frac{\eps}{3})|V(G_x)|$ for all $x\in X$.
By Theorem~\ref{thm:rrsmatching}, $G_x$ then has a perfect matching $M_x$ for all $x\in X$.
Let $A_2 := \bigcup_{x\in X}M_x$.
Because sets in $\cQ$ are pairwise disjoint, for $x,y \in X$, if $f(x)$ and $f(y)$ are disjoint, then $M_x$ and $M_{y}$ are vertex-disjoint.

Let $\cP_2$ be a refinement of $\cQ$ into $k$-sets, and let $\cP := \cP_1 \cup \cP_2$.
Let $H_0$ be a $\cP_1$-compatible $k$-graph with $V(H_0) \subseteq U_1$, and suppose that $A_1+H_0$ admits a decomposition $\cM$ into $\cP_1$-compatible matchings.
Then $A_1$ is $\cP_1$-compatible, so for $P \in \cP_1$, let $d_P := d_{A_1}(v)$, where $v\in P$.
For each $e\in \binom{U_1}{k}$, let $\cM_e := \{M \in \cM: e\in M\}$ be the matchings in $\cM$ that contain $e$.

Let $f_0:X\cup E(H_0) \to \binom{U_1}{k}$ be the function with $f_0(x)=f(x)$ and $f(e)=e$ for all $x\in X$ and $e\in E(H_0)$ (so that $f_0$ represents the $k$-multigraph $A_1+H_0$ by combining the ``fake edge" map $f$ with the identity map on $E(H_0)$).
Let $g:X\cup E(H_0) \to \cM$ be a function such that for all $e\in \binom{U_1}{k}$, $g$ maps $f_0^{-1}(e)$ bijectively onto $\cM_e$, which is possible because the number $|\cM_e|$ of matchings in $\cM$ that contain $e$ is the multiplicity $|f_0^{-1}(e)|$ of $e$ in $A_1+H_0$.
Thus, each fake edge $x\in X$ has two associated matchings: the matching $M_x$ in $G_2$, and the matching $g(x)$ in $A_1+H_0$.

For each $M\in \cM$, let
\[M' := (M\cap g^{-1}(M)) \cup \bigcup_{x\in X \cap g^{-1}(M)}M_x.\]
$M \cap g^{-1}(M)$ represents the part of $M$ that comes from $H_0$, and $\bigcup_{x \in X\cap g^{-1}(M)}M_x$ is the union of the matchings that replace the fake edges in $M$.

We claim that $\cM' := \{M':M\in \cM\}$ is a decomposition of $H_0 \cup A_2$ into $\cP$-compatible matchings.
Indeed, the $M_x$ are pairwise edge-disjoint and each edge-disjoint from $H_0$, since each edge of $M_x$ has at least one vertex in $Q_x$, which does not intersect $H_0$ or any $M_y$ with $y \neq x$.
Every edge $e \in E(H_0)$ is in the matching $g(e)' \in \cM'$ and no other matching in $\cM'$, and every edge in $A_2$ is in exactly one $M_x$ for some $x\in X$ and thus in exactly one $M' \in \cM'$.
Thus, $\cM'$ is a decomposition of $H_0 \cup A_2$ into subgraphs.

Let $M \in \cM$, and observe that to show that $M'$ is a $\cP$-compatible matching, it suffices to show that the multiset $L:=V(M\cap g^{-1}(M)) + \sum_{x\in X\cap g^{-1}(M)}V(M_x)$ is a $\cP$-compatible set (with multiplicities $\le 1$), since the multiplicities in $L$ are precisely the degrees of vertices in $M'$.
For each $e\in M\setminus g^{-1}(M)$, there is a unique $x\in f^{-1}(e)$ such that $g(x) = M$ by the definition of $g$, so $V(M\setminus g^{-1}(M)) = \sum_{x\in X \cap g^{-1}(M)}f(x)$ as multisets.
Thus, since the $Q_x$ are pairwise disjoint and disjoint from $U_1 \supseteq V(M)$, we have
\[ L = V(M\cap g^{-1}(M))+\sum_{x\in X\cap g^{-1}(M)}(Q_x \cup f(x)) = V(M) \cup \bigcup_{x\in X\cap g^{-1}(M)}Q_x, \]
which is a $\cP$-compatible set, since $M$ is $\cP$-compatible and sets in $\cQ$ are $\cP$-compatible.
Therefore, $\cM'$ is indeed a decomposition of $H_0\cup A_2$ into $\cP$-compatible matchings, as desired.
\end{proof}

Finally, the following lemma details the third step of the construction of the absorber, which completes the compatible matchings from the previous step to perfect matchings (on the vertex set of the absorber).

\begin{lemma} \label{lem:absorbingperfect}
    Let $k \ge 3$.
    Let $\eps>0$ and $m_2,m_3 \in k\bN$ such that $1/m_3 \ll \eps$ and $m_2 \le m_3^{1/3}$.
    Let $U_2 \subseteq U_3$ be sets with $|U_2| = m_2$ and $|U_3| = m_3$, let $G_3$ be a $k$-graph with $V(G_3) = U_3$ and $\delta_{k-1}(G_3) \ge (\frac12 + \eps)m_3$, and let $\cP$ be a partition of $U_2$ into $k$-sets.
    For each $P\in \cP$, let $0\le d_P \le \binom{m_2-1}{k-1}$.
    Then there exists a spanning subgraph $A_3$ of $G_3$ such that $U_2$ is an independent set in $A_3$ and for any $\cP$-decomposable $k$-graph $H_2$ with $V(H_2) \subseteq U_2$ and $d_{H_2}(v)=d_P$ for all $v\in P \in \cP$, $A_3 \cup H_2$ admits a $1$-factorization.
\end{lemma}
\begin{proof}
Choose a new constant $\eps_1$ such that $1/m_3 \ll \eps_1 \ll \eps$.
Let $\cQ = \{Q_P:P \in \cP\}$ be a random $k$-equitable $\frac{m_2}{k}$-partition of $U_3\setminus U_2$ chosen such that each $Q_P$ is a uniformly random subset of $U_3\setminus U_2$ of a fixed size.
Since $m_2/m_3 \le m_3^{-2/3} \ll \eps/4$ and $|Q_P| = \frac{k(m_3-m_2)}{m_2}\pm k = \Omega(m_3^{2/3})$ for all $P\in \cP$, by Lemma~\ref{lem:randomset} with $Z=U_2$ and a union bound, the probability that $d_{G_3,Q_P}(S) \le (\frac12+\frac\eps2)|Q_P|$ for some $P \in \cP$ and $S\in \binom{m_3}{k-1}$ is at most
\[ \frac{m_2}{k} \binom{m_3}{k-1}e^{-\Omega\left(\eps^2m_3^{2/3}\right)} <1, \]
so fix a choice of $\cQ$ such that $\delta_{k-1}(G_3[Q_P\cup P]) \ge (\frac12 +\frac\eps2)|Q_P|$ for all $P \in \cP$.
In particular, for $P\in \cP$, letting $G_P$ be $G_3[Q_P \cup P]$ with the edge $P$ removed (if it exists), we have $\delta_{k-1}(G_P) \ge (\frac12 +\frac{\eps}{3})|V(G_P)|$.

Let $d_3 := \binom{m_2}{k}$.
For each $P\in \cP$, since $|V(G_P)| = \Omega(m_3^{2/3})$, we have $d_3 \le m_3^{k/3} \ll \eps_1 |V(G_P)|^{k-1}$ and $|P| = k \ll \eps_1 |V(G_P)|$.
Thus, by Lemma~\ref{lem:sparsematchings} with $\delta= \eps/3$ and $\delta_1 = \eps_1$, there exist pairwise edge-disjoint matchings $M_{P,1},\ldots,M_{P,d_3}$ in $G_3$ such that $V(M_{P,i}) = Q_P$ for $i\le d_P$ and $V(M_{P,i}) = Q_P \cup P$ for $i > d_P$.
Let $A_3$ be the spanning subgraph of $G_3$ with $E(A_3) = \bigcup_{P \in \cP}\bigcup_{i=1}^{d_3}M_{P,i}$, and note that because $U_2 \cap V(G_P) = P$ is an independent set in each $G_P$, $U_2$ is an independent set in $A_3$.

\begin{figure}
    \centering
    \includegraphics{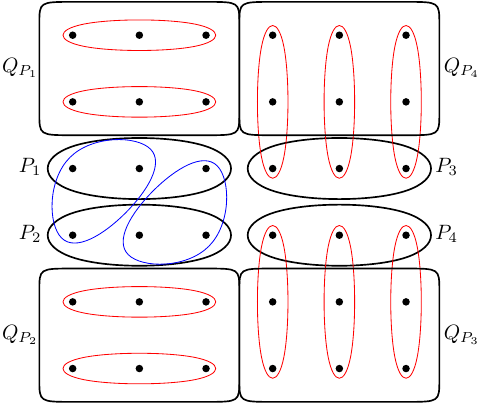}
    \caption{Absorbing a $\cP$-compatible matching (in blue) using matchings on $Q_P$ and on $Q_P\cup P$ for $P\in \cP$ (in red) into a perfect matching in Lemma~\ref{lem:absorbingperfect}, where $\cP = \{P_1,P_2,P_3,P_4\}$.}
    \label{fig:absorbingperfect}
\end{figure}

Let $H_2$ be a $k$-graph with $V(H_2) \subseteq U_2$ such that $d_{H_2}(v)=d_P$ for all $v\in P \in \cP$ and $H_2$ admits a decomposition $\cM = \{M_1,\ldots,M_t\}$ into (non-empty) $\cP$-compatible matchings.
Since $t \le d_3$, let $M_{t+1} = \cdots = M_{d_3} := \emptyset$.
For each $P\in \cP$, exactly $d_P$ matchings in $\cM$ cover $P$, so let $\sigma_P:[d_3]\to[d_3]$ be a bijection such that $M_{\sigma_P(1)},\ldots,M_{\sigma_P(d_P)}$ are the matchings in $\cM$ that cover $P$.
For each $i\in [d_3]$, let
\[ M_i' := M_i \cup \bigcup_{P \in \cP}M_{P,\sigma_P^{-1}(i)}. \]
We claim that $\cM' := \{M_1',\ldots,M_{d_3}'\}$ is a decomposition of $A_3 \cup H_2$ into perfect matchings.
Indeed, consider the following equality of vertex multisets:
\[ V(M_i) + \sum_{P \in \cP}V(M_{P,\sigma_P^{-1}(i)}) = \sum_{\substack{P \in \cP \\ \sigma_P^{-1}(i) \le d_P}}P+\sum_{\substack{P \in \cP \\ \sigma_P^{-1}(i) \le d_P}}Q_P+\sum_{\substack{P \in \cP \\ \sigma_P^{-1}(i) > d_P}}(Q_P+P) = V(\cP)+ V(\cQ) = U_3. \]
This holds since $M_i$ is $\cP$-compatible, and $P \in \cP$ has $P \subseteq M_i$ if and only if $i \in \{\sigma_P(1),\ldots,\sigma_P(d_P)\}$, i.e., $\sigma_P^{-1}(i) \le d_P$.
Since the multiplicities in this multiset are precisely the degrees of vertices in $M_i'$, $M_i'$ is $1$-regular, i.e., a perfect matching.
Also, $\sum_{i=1}^{d_3}M_i' = E(A_3 \cup H_2)$ as edge multisets by the definitions of $A_3$ and $\cM$.
Therefore, $\cM'$ is a decomposition of $A_3 \cup H_2$ into perfect matchings, as desired.
\end{proof}

\subsection{Proof of main absorption lemma}
\begin{proof}[Proof of Lemma~\ref{lem:absorbing}]
Assume $\cU$ contains at least one $m$-set (else redefine $m:= m+k$), and let $N:= |\cU|$ and $\cU = \{U_{1,0},\ldots,U_{N,0}\}$.
Note that $\frac{n}{m+k} < N \le \frac{n}{m}$.
Let $m_1 := k \lflr\frac{n^{0.05}}{k}\rflr$, $m_2 := m_1^2$, and $m_3 := m_1^6$, so that $(1-o(1))n^{0.1} \le m_2 \le n^{0.1}$ and $(1-o(1))n^{0.3}\le m_3 \le n^{0.3}$.
Let $\delta_1,\delta_2$ be new constants with $\delta_1 \ll \delta_2 \ll \gamma$.
If $d_0>\binom{m-1}{k-1}$, then the conclusion holds vacuously for empty $A$, so assume $d_0 \le \binom{m-1}{k-1}$.

For $i \in [N]$, let $m_{i,0} := |U_{i,0}| \in \{m,m+k\}$, and let $U_{i,1},U_{i,2},U_{i,3}$ be uniformly random subsets of $V(G)$ conditioned on $U_{i,0} \subseteq U_{i,1} \subseteq U_{i,2} \subseteq U_{i,3}$ and $|U_{i,j}| = m_j$ for all $j \in [3]$, chosen independently for each $i$.
For $i \in [N]$ and $j\in \{2,3\}$, let $B_{i,j}$ be the event that $\delta_{k-1}(G[U_{i,j}]) < (\frac12+\frac{\delta}{3})m_j$.
By~Lemma~\ref{lem:randomset} with $\eps = \delta$, since $m_j/n = O(n^{-0.7}) \le \delta/4$ and $U_{i,j}\setminus U_{i,0}$ is a uniformly random $(m_j-m)$- or $(m_j-(m+k))$-subset of $V(G)\setminus U_{i,0}$, we have that $\prob{B_{i,j}} \le e^{-\Omega(\delta^2m_j)}$.
In particular, 
\begin{equation} \label{eq:badevents}
    \prob{\bigcup_{i=1}^N(B_{i,2}\cup B_{i,3})} \le 3N e^{-\Omega(\delta^2n^{0.1})} = o(1).
\end{equation}

For $S\subseteq V(G)$, let $X_S := \sum_{i=1}^N\IND_{S \subseteq U_{i,3}}$, and for $v\in V(G)$, let $X_v := X_{\{v\}}$.

\begin{claim} \label{clm:setinfewrandomsets}
    W.h.p. for all $S\in \binom{V(G)}{k-1}$, $X_S \le \log^2{n}$.
\end{claim}
\begin{proof}[Proof of Claim]
Let $S \in \binom{V(G)}{k-1}$, and let $I_S := \{i \in [N]:S\cap U_{i,0} = \emptyset\}$.
Let $Y_S := \sum_{i\in I_S}\IND_{S\subseteq U_{i,3}}$, and note that since the $U_{i,0}$ are pairwise disjoint, $S$ intersects at most $k-1$ of them, so $|[N]\setminus I_S| \le k-1$, and thus $X_S \le Y_S + k-1$.
Since $k\ge 3$ and $U_{i,3}\setminus U_{i,0}$ is a uniformly random $(m_3-m_{i,0})$-subset of $V(G)\setminus U_{i,0}$,
\begin{align*}
    \expn{Y_S} &\le \sum_{i \in I_S}\prob{S \subseteq U_{i,3}} = |I_S| \cdot \left(\frac{(m_3-m_{i,0})(m_3-m_{i,0}-1)\cdots (m_3-m_{i,0}-(k-2))}{(n-m_{i,0})(n-m_{i,0}-1)\cdots (n-m_{i,0}-(k-2))}\right) \\
    &= O\left(\frac{n}{m}\left(\frac{m_3}{n}\right)^{k-1}\right) = O\left(\frac{n^{1.7-0.7k}}{m}\right) = O\left(n^{-0.4}\right).
\end{align*} 
Thus, by Theorem~\ref{thm:chernoffhyper}, since $\IND_{S\subseteq U_{i,3}}$ are independent for $i \in I_S$,
\[ \prob{X_S > \log^2n}\le \prob{Y_S > \frac12 \log^2{n}} \le e^{-\Omega(\log^2{n})}. \]
Therefore, the probability that $X_S> \log^2{n}$ for some $S \in \binom{V(G)}{k-1}$ is at most
\[ n^{k-1}e^{-\Omega(\log^2{n})} = n^{k-1-\Omega(\log{n})} = o(1), \]
proving the claim.
\end{proof}

\begin{claim} \label{clm:absorberscodegree}
    W.h.p. for all $i \in [N]$, $\Delta_{k-1}(\bigcup_{i'=1}^{i-1}G[U_{i',3} \cap U_{i,3}]) \le 2m\log^2n$.
\end{claim}
\begin{proof}[Proof of Claim]
Let $i \in [N]$ and $S \in \binom{U_{i,3}}{k-1}$.
By Claim~\ref{clm:setinfewrandomsets}, w.h.p. $X_S \le \log^2n$, so let $I \subseteq [i-1]$ with $|I| \le \log^2n$ and condition on $\{i' \in [i-1]:S \subseteq U_{i',3}\} = I$.
Then
\[d_{\bigcup_{i'=1}^{i-1}G[U_{i',3} \cap U_{i,3}]}(S) \le \sum_{i' \in I}|U_{i',3}\cap (U_{i,3} \setminus S)| \le (m+k)\log^2n+\sum_{i' \in I}|(U_{i',3}\setminus U_{i',0})\cap (U_{i,3}\setminus S)|. \]
Conditioned on a fixed choice for $U_{i,3}$, the variables $Z_{S,i'} := |(U_{i',3}\setminus U_{i',0})\cap (U_{i,3}\setminus S)|$ for $i' \in I$ are independent hypergeometric random variables with $\expn{Z_{S,i'}} = O(\frac{m_3^2}{n}) = O(n^{-0.4})$, since $U_{i',3}\setminus (U_{i',0} \cup S)$ is a uniformly random $(m_3-|U_{i',0}\cup S|)$-subset of $V(G)\setminus (U_{i',0} \cup S)$ and $U_{i,3}\setminus (U_{i',0} \cup S)$ is a fixed $O(m_3)$-subset of $V(G)\setminus U_{i',0}$.
Thus, by Theorem~\ref{thm:chernoffhyper}, letting $Z_S := \sum_{i' \in I}Z_{S,i'}$,
\begin{align*}
    \prob{d_{\bigcup_{i'<i}G[U_{i',3} \cap U_{i,3}]}(S)>2m\log^2n} &\le\prob{(m+k)\log^2n+Z_S > 2m\log^2n} \\
    &\le \prob{Z_S \ge \log^2n} \le e^{-\Omega(\log^2n)}.
\end{align*}
Therefore, the probability that $\Delta_{k-1}(\bigcup_{i'=1}^{i-1}G[U_{i',3} \cap U_{i,3}]) > 2m\log^2n$ for some $i \in [N]$ is at most
\[ N\binom{m_3}{k-1}e^{-\Omega(\log^2n)} = e^{O(\log{n})-\Omega(\log^2n)} = o(1), \]
proving the claim.
\end{proof}

For each $i \in [N]$, let $\cP_{i,1}$ be a partition of $U_{i,1}$ into $k$-sets such that $U_{i,0}$ is $\cP_{i,1}$-compatible.
We will now put Lemmas~\ref{lem:absorbingmulti}, \ref{lem:absorbingcompatible}, and \ref{lem:absorbingperfect} together to construct our desired absorber.
\begin{claim} \label{clm:absorber}
    W.h.p. there exist pairwise edge-disjoint $k$-graphs $A_1,\ldots,A_N$ such that for each $i\in [N]$, we have that $A_i \subseteq G[U_{i,3}] - G[\cU]$ and for any $d_0$-regular $k$-graph $H_{i,0}$ with $V(H_{i,0}) = U_{i,0}$, $A_i \cup H_{i,0}$ admits a decomposition into $d_i$ perfect matchings for some $d_i \in \bN$.
\end{claim}
\begin{proof}[Proof of Claim]
By~\eqref{eq:badevents}, w.h.p. $\overline{B_{i,j}}$ holds for all $i \in [N]$ and $j\in \{2,3\}$, i.e., $\delta_{k-1}(G[U_{i,j}]) \ge (\frac12+\frac\delta3)m_j$; assume this occurs.
Also assume that the event in Claim~\ref{clm:absorberscodegree} occurs.
We will construct $A_1,\ldots,A_N$ inductively.
Let $i \in [N]$ and suppose that $A_1,\ldots,A_{i-1}$ have already been defined.
Let $H_{i,0}$ be a $d_0$-regular $k$-graph with $V(H_{i,0}) = U_{i,1}$.

By Lemma~\ref{lem:absorbingmulti} with $\cP = \cP_{i,1}$, since $m_1 = \Theta(n^{0.05}) \gg m_{i,0}+\binom{m_{i,0}}{k}k(k-1)$, there exists a $k$-multigraph $A_{i,1}$, not depending on the choice of $H_{i,0}$, with $V(A_{i,1}) = U_{i,1}$ and $|E(A_{i,1})| = O(m^{k^2})$ such that $A_{i,1}+H_{i,0}$ is $\cP_{i,1}$-decomposable.

By Claim~\ref{clm:absorberscodegree}, we can bound the $(k-1)$-degree of the previously defined $A_{i'}$ as follows:
\[ \Delta_{k-1}\left(\bigcup_{i'=1}^{i-1}A_{i'}[U_{i,2}]\right) \le \Delta_{k-1}\left(\bigcup_{i'=1}^{i-1}A_{i'}[U_{i,3}]\right) \le \Delta_{k-1}\left(\bigcup_{i'=1}^{i-1}G[U_{i',3}\cap U_{i,3}]\right) \le 2m\log^2n. \]
We also have $\Delta_{k-1}(G[\cU]) \le m+k$.
So, for $j \in \{2,3\}$, because $\overline{B_{i,j}}$ holds, the $k$-graph $G_j := G[U_{i,j}]-\bigcup_{i'=1}^{i-1}A_{i'} - G[\cU]$ has 
\[\delta_{k-1}(G_j) \ge \left(\frac12+\frac{\delta}{3}\right)m_j-2m\log^2n - (m+k) \ge \left(\frac12+\frac{\delta}{4}\right)m_j. \]

By Lemma~\ref{lem:absorbingcompatible} with $\eps = \delta/4$ applied to $G_2$, since $|E(A_{i,1})| = O(m^{k^2}) \ll \sqrt{m_2}$, there exists a partition $\cP_{i,2}$ of $U_{i,2}\setminus U_{i,1}$ into $k$-sets and a subgraph $A_{i,2} \subseteq G_2$ with $V(A_{i,2}) = U_{i,2}$, not depending on the choice of $H_{i,0}$, such that $U_{i,1}$ is an independent set in $A_{i,2}$ and, since $A_{i,1}+H_{i,0}$ is $\cP_{i,1}$-decomposable, $A_{i,2} \cup H_{i,0}$ is $(\cP_{i,1} \cup \cP_{i,2})$-decomposable.

By Lemma~\ref{lem:absorbingperfect} applied to $G_3$ with $\eps = \delta/4$, $\cP = \cP_{i,1}\cup \cP_{i,2}$ and $d_P = d_{A_{i,2}}(v)+d_0\IND_{P \subseteq U_{i,0}}$ for all $v\in P \in \cP$, there exists a spanning subgraph $A_{i,3} \subseteq G$ with $V(A_{i,3}) = U_{i,3}$, not depending on the choice of $H_{i,0}$, such that, since $A_{i,2}\cup H_{i,0}$ is $\cP$-decomposable, $A_{i,3} \cup (A_{i,2} \cup H_{i,0})$ admits a $1$-factorization.
Therefore, the claim follows by setting $A_i := A_{i,2} \cup A_{i,3}$.
\end{proof}

By Claim~\ref{clm:absorber}, w.h.p. $A' := \bigcup_{i=1}^N A_i$ exists.
Observe that by Claim~\ref{clm:setinfewrandomsets}, w.h.p. for all $S \in \binom{V(G)}{k-1}$,
\[ d_{A'}(S) =  \sum_{i=1}^N\IND_{S \subseteq U_{i,3}}\cdot d_{A_i}(S) \le m_3X_S  \le n^{0.3}\log^2{n}, \]
so 
\[\delta_{k-1}(G-A'-G[\cU]) \ge \left(\frac12+\delta\right)n-n^{0.3}\log^2n-(m+k) \ge \left(\frac12+\frac{\delta}{2}\right)n. \]
By \eqref{eq:badevents}, observe that w.h.p. for all $i \in [N]$, $|U_{i,3}| = m_3 = O(n^{0.3}) \ll \delta_1 n$, $k$ divides $n-|U_{i,3}|$, and $d_i \le \binom{m_3-1}{k-1}$, which implies that $\ell := \sum_{i=1}^N d_i \le \frac{m_3^{k-1}n}{m} \le n^{0.3k+0.7} \ll \delta_1n^{k-1}$.
Thus, by Lemma~\ref{lem:sparsematchings} applied to $G-A'-G[\cU]$ with $\delta$ replaced by $\delta/2$, w.h.p. there exist pairwise edge-disjoint matchings $M_{i,j}$ in $G-A'-G[\cU]$ for $i\in [N]$ and $j \in [d_i]$ such that $V(M_{i,j}) = V(G)\setminus U_{i,3}$ and $\Delta_{k-1}(\bigcup_{i=1}^n\bigcup_{j=1}^{d_i}M_{i,j}) \le \delta_2 n$.

Then, letting $A := A' \cup \bigcup_{i=1}^N\bigcup_{j=1}^{d_i}M_{i,j}$, we have that $A$ is a spanning subgraph of $G-G[\cU]$ with
\[ \Delta_{k-1}(A) \le n^{0.3}\log^2{n}+\delta_2n \le \gamma n. \]
Because $A$, as defined above, exists w.h.p., and $n$ is sufficiently large, fix such a choice of $A$.

It remains to show that $A$ has the desired absorbing property.
Let $K$ be a $d_0$-regular spanning subgraph of $K_n^{(k)}[\cU]$.
For each $i\in [N]$, $K[U_{i,0}]$ is $d_0$-regular, so by Claim~\ref{clm:absorber}, $A_i \cup K[U_{i,0}]$ admits a decomposition into matchings $M_{i,1}',\ldots,M_{i,d_i}'$.
For each $i \in [N]$ and $j\in [d_i]$, $M_{i,j} \cup M_{i,j}'$ is a perfect matching in $G$, and these matchings are pairwise edge-disjoint and contain all edges of $A \cup K$, so $A\cup K$ admits a $1$-factorization.
Therefore, $A$ satisfies the conditions that $\Delta_{k-1}(A)\le \gamma n$ and $A\cup K$ admits a $1$-factorization for any $d_0$-regular spanning subgraph $K$ of $K_n^{(k)}[\cU]$, completing the proof of Lemma~\ref{lem:absorbing}.
\end{proof}

\section{Proof of main theorem} \label{sec:main}
\begin{proof}[Proof of Theorem~\ref{thm:main}]
Let $C, \gamma, L$ be new constants such that $0 < 1/n \ll 1/C \ll \gamma \ll 1/L \ll \delta$, and let $\alpha := 1/L$.
Let $r$ be the vertex degree of $G$ and $\beta := r/n^{k-1}$.
By Lemma~\ref{lem:vortex}, there exists a $(\delta,\beta,L,C)$-vortex $\cP_0,\ldots,\cP_{t+1}$ for $G$.
Observe that
\[\delta_{k-1}(G-G[\cP_1]) \ge \delta_{k-1}(G)- \max_{V\in\cP_1}|V| \ge \left(\frac12+\delta-\alpha\right)n-k \ge \left(\frac12+ \frac{\delta}{2}\right)n. \]
By Corollary~\ref{cor:partitionsequence}, $|U| \in \left\{k\ff{\alpha^tn}{k},k\cf{\alpha^tn}{k}\right\}$ for all $U \in \cP_t$.
Let $d_0 := \lcl \frac12 \beta (\alpha^t n)^{k-1} \rcl$.

Apply Lemma~\ref{lem:absorbing} with $m=k\lflr \frac{\alpha^tn}{k}\rflr \le L^2C$, $\cU = \cP_t$, $\delta/2$ in place of $\delta$, and $G-G[\cP_1]$ in place of $G$ to obtain a spanning subgraph $A$ of $G-G[\cP_1]$ such that $\Delta_{k-1}(A) \le \gamma n$ and for any $d_0$-regular spanning subgraph $K$ of $K_n^{(k)}[\cP_t]$, $A\cup K$ admits a decomposition into $d_1$ perfect matchings for some $d_1\in \bN$.
Since at least one such $K$ exists (by taking $d_0$ matchings from an application of Baranyai's theorem to the complete graph on $U$ for each $U \in \cP_t$ and taking their union across all $U \in \cP_t$), it follows that $A$ is $(d_1-d_0)$-regular and $d_1$ does not depend on the choice of $K$.

Observe that then $\cP_1$ is a $(\delta/2,\beta,L,1)$-partition for $G-A$, since
\[ d_{G-A,U}(S) \ge d_{G,U}(S)-\gamma n \ge \left(\frac12+\frac{\delta}{2} - |U|^{-1/3}\right)|U| \]
for $U \in \cP_1$ and $S \in \binom{V(G)}{k-1}$, and $(G-A)[\cP_1] = G[\cP_1]$.
Similarly, $(G-A)[\cP_s] = G[\cP_s]$ for all $s\in [t+1]$, so $\cP_0,\ldots,\cP_{t+1}$ is a $(\delta/2, \beta,L,C)$-vortex for $G-A$.
Thus, by Lemma~\ref{lem:iterative} with $\delta$ replaced by $\delta/2$, there exists a $k$-graph $H$ such that $G-A-G[\cP_t] \subseteq H \subseteq G-A-G[\cP_{t+1}]$, $\Delta_{k-1}(H[\cP_t]) \le \alpha^{t+2}n$, and $H$ admits a decomposition into $\cS$-matchings, where $\cS := \bigcup_{s=0}^t \cP_s$.
Because each of these $\cS$-matchings is $\cP_t$-compatible, $H$ is $\cP_t$-compatible, i.e., $d_H(u)=d_H(v)$ whenever $u,v \in U \in \cP_t$.

Let $F := G-A-H$, and note that because $G,A,H$ are $\cP_t$-compatible, so is $F$, and $F\subseteq G[\cP_t]$, meaning that for each $U\in \cP_t$, $F[U]$ is $d_U$-regular for some $d_U$.
For $v\in U \in \cP_t$,
\[ d_U = d_F(v) = d_{G[U]}(v)-d_{H[U]}(v) = \left(1\pm |U|^{-1/3}\right)\beta |U|^{k-1}- O\left(|U|^{k-2}\alpha^{t+2}n\right) = \left(1\pm O\left(\alpha^2\right)\right)\beta |U|^{k-1}. \]
For each $U\in \cP_t$,
\[\delta_{k-1}(F[U]) \ge \delta_{k-1}(G[U])-\Delta_{k-1}(H[U]) \ge \left(\frac12+\delta - |U|^{-1/3}\right)|U|-\alpha^{t+2}n \ge \left(\frac12+\gamma\right)|U|, \]
so by Theorem~\ref{thm:approxdecomp}, $F[U]$ contains $(1-\gamma)d_U = (1\pm O(\alpha^2))\beta |U|^{k-1}$ edge-disjoint perfect matchings.
Let $F_U$ be the union of $d_U-d_0$ of these matchings (since $d_U-d_0 < (1-\gamma)d_U$), and let $F' := \bigcup_{U\in \cP_t}F_U$.
Then $F'$ admits a decomposition into $\cS$-matchings.

Now, $K:= F-F'$ is a $d_0$-regular subgraph of $G[\cP_t]$, so $A \cup K$ admits a $1$-factorization.
Since $H$ and $F'$ each admit decompositions into $\cS$-matchings, $G = A \cup H \cup F' \cup K$ admits a decomposition into $\cS$-matchings.
Therefore, by Lemma~\ref{lem:S-matching}, $G$ admits a $1$-factorization, completing the proof of Theorem~\ref{thm:main}.
%
%
\end{proof}

\section{Lower bound construction} \label{sec:lowerbound}
The following result adapts the constructions by R\"odl, Ruci\'nski, and Szemer\'edi~\cite{RRS:09} (see also \cite{KO:06}) of $1$-factor-free $k$-graphs with large minimum $(k-1)$-degree.
While the constructions in~\cite{RRS:09} are not necessarily vertex-regular, we adjust the constructions so that they are vertex-regular and still admit no $1$-factorization.

\begin{prop} \label{prop:lowerbound}
    For all $k\ge 3$, there exist infinitely many $n\in k\bN$ for which there exists a vertex-regular $k$-graph $G$ on $n$ vertices with $\delta_{k-1}(G) \ge \frac12 n - O_k(\sqrt{n})$ that admits no $1$-factorization.
\end{prop}
\begin{proof}
First, consider the case in which $k$ is even.
Let $n \in \bN$ with $n \equiv k \pmod{2k}$, and let $A,B$ be disjoint $(n/2)$-sets.
Then $|A|=n/2 \equiv k/2 \pmod{2}$.
Let $G$ be the $k$-graph on the vertex set $A\cup B$ whose edges consist of all $e \in \binom{A\cup B}{k}$ with $|e\cap A| \equiv k/2+1\pmod{2}$.
If $k/2$ is even, then any perfect matching in $G$ would have $n/k$ edges, and $n/k$ is odd, so it would cover an odd number of vertices in $A$, but $|A|$ is even.
If $k/2$ is odd, then every matching in $G$ would cover an even number of vertices in $A$, but $|A|$ is odd.
Thus, $G$ has no perfect matching.
Also, $G$ is vertex-regular, since swapping the roles of the parts $A$ and $B$ results in the same $k$-graph, and vertices in the same part have the same degree.
Further, $\delta_{k-1}(G) \ge \frac12 n -k$, since each $S \in \binom{A\cup B}{k-1}$ forms an edge with all vertices in either $A \setminus S$ or $B\setminus S$.
Therefore, we are done in this case.

Next, consider the case in which $k$ is odd.
Let $n \in k\bN$, and let $C$ be a new constant with $1/n \ll 1/C \ll 1/k$.
Let $A$ and $B$ be disjoint sets, and let $H$ be the $k$-graph on $A\cup B$ with $E(H) = \{e:|e\cap A| \text{ is even}\}$.
Then all vertices in $A$ have the same degree $d_A$, and all vertices in $B$ have the same degree $d_B$.
Let $a := |A|$, $b := |B|$, and $\ell := k-1$.
Denote the coefficient of $x^t$ in a polynomial $f$ by $[x^t]f$.
Assuming $a\le b$, applying Pascal's formula and reindexing,
\begingroup
\allowdisplaybreaks
\begin{align}
    f(a,b) &:= d_B-d_A = \sum_{\substack{j \le \ell \\ j \text{ even}}}\binom{a}{j}\binom{b-1}{\ell-j}-\sum_{\substack{j \le \ell \\ j \text{ odd}}}\binom{a-1}{j}\binom{b}{\ell-j} \nonumber \\
    &= \sum_{\substack{j \le \ell \\ j \text{ even}}}\left(\binom{a-1}{j}+\binom{a-1}{j-1}\right)\binom{b-1}{\ell-j}-\sum_{\substack{j \le \ell \\ j \text{ odd}}}\binom{a-1}{j}\left(\binom{b-1}{\ell-j}+\binom{b-1}{\ell-j-1}\right) \nonumber \\
    &= \sum_{j=0}^{\ell} (-1)^j\binom{a-1}{j}\binom{b-1}{\ell-j}+\sum_{\substack{2\le j \le \ell \\ j \text{ even}}}\binom{a-1}{j-1}\binom{b-1}{\ell-j}-\sum_{\substack{1\le j \le \ell-1 \\ j \text{ odd}}}\binom{a-1}{j}\binom{b-1}{\ell-j-1} \nonumber \\
    &= \sum_{j=0}^\ell(-1)^j\binom{a-1}{j}\binom{b-1}{\ell-j} = [x^\ell] (1-x)^{a-1}(1+x)^{b-1}  \nonumber \\
    &= [x^\ell]\left(1-x^2\right)^{a-1}(1+x)^{b-a} = \sum_{j=0}^{\lflr \ell/2\rflr}(-1)^j\binom{a-1}{j}\binom{b-a}{\ell-2j}. \label{eq:degreedifference}
\end{align}
\endgroup
\begin{claim} \label{clm:constructionsizes1}
    Suppose $k \equiv 1 \pmod{4}$ and $0\le j < 2k$.
    Then there exist $a,b$ with $a+b=n$, $a \equiv j \pmod{2k}$, $0\le b-a \le 4k$, and $f(a,b) = \Theta_k\left(n^{\ell/2}\right)$.
\end{claim}
\begin{proof}[Proof of Claim]
Let $n/2-2k < a\le n/2$ with $a \equiv j \pmod{2k}$ and $b=n-a$.
Then $0\le b-a \le 4k$, so by~\eqref{eq:degreedifference},
\[ f(a,b) = (-1)^{\ell/2}\binom{a-1}{\ell/2} \pm O_k\left(a^{\ell/2-1}\right) = \Theta_k\left(n^{\ell/2}\right). \qedhere \]
\end{proof}

\begin{claim} \label{clm:constructionsizes2}
    Suppose $k \equiv 3\pmod{4}$ and $0\le j < 2k$.
    Then there exist $a,b$ with $a+b=n$, $a \equiv j \pmod{2k}$, $0 \le b-a \le Cn^{1/2}$, and $0\le f(a,b) = O_C(n^{(\ell-1)/2})$.
\end{claim}
\begin{proof}[Proof of Claim]
Let $h:= \frac12(n-C\sqrt{n})$, and let $a \in [h,n/2] =: I$ be maximal such that $a \equiv j \pmod{2k}$ and, setting $b=n-a$, we have $d_B-d_A \ge 0$.
This is possible because for $h \le a <h+2k$, $b-a = \Theta(C\sqrt{n})$, so by~\eqref{eq:degreedifference},
\[ f(a,b) = \binom{b-a}{\ell} \pm O_k\left(C^{\ell-2}n^{\ell/2}\right) = \Theta_k\left(C^\ell n^{\ell/2}\right) \ge 0. \]
Now, $b-a \le C\sqrt{n}$ by definition, and because $f(x,n-x)$ can be treated as a polynomial in $x$ with a zero in $I$ (since $f(h,n-h)\ge 0$ and $f(n/2,n/2) = (-1)^{\ell/2} \binom{n/2-1}{\ell/2} = -\binom{n/2-1}{\ell/2} \le 0$) and $a$ is within distance at most $2k$ from such a zero,
\[ 0\le f(a,b) \le 2k\cdot \sup_{x \in I}\left|\frac{d}{d x}f(x,n-x)\right| = O_k\left(C^\ell n^{(\ell-1)/2}\right). \qedhere \]
\end{proof}


Choose $a,b$ as in Claim~\ref{clm:constructionsizes1} or Claim~\ref{clm:constructionsizes2} with $a\equiv k \pmod{2k}$.
Observe that $H$ has no perfect matchings, since any matching in $H$ has even intersection with $A$, but $a$ is odd.
Also note that $A$ is an independent set in $H$, since edges contained entirely in $A$ do not have an even intersection with $A$.

Because $k \mid a$ and $0\le f(a,b) = O(n^{\ell/2})$, by Baranyai's theorem, the complete $k$-graph on the vertex set $A$ has a decomposition into $\binom{a-1}{\ell} \ge f(a,b)$ perfect matchings.
By adding $f(a,b)$ of these to $H$, we obtain a vertex-regular $k$-graph $G$ by~\eqref{eq:degreedifference}.
Further, since $H$ has no perfect matching, each perfect matching in $G$ must contain one of these added edges, so $G$ has at most $\frac{a}{k}f(a,b)$ edge-disjoint perfect matchings, which together contain at most $\frac{na}{k^2}f(a,b)$ edges.
If $k\equiv 1 \pmod{4}$, then by Claim~\ref{clm:constructionsizes1}, this is at most $O_k\left(n^2n^{\ell/2}\right) = O_k\left(n^{(k+3)/2}\right) = o(n^k)$ edges, and if $k\equiv 3\pmod{4}$, then by Claim~\ref{clm:constructionsizes2}, this is at most $O_C\left(n^2n^{(\ell-1)/2}\right) = O_C\left(n^{(k+2)/2}\right) = o(n^k)$ edges.
But $G$ has $\Omega(n^k)$ edges, so $G$ has no $1$-factorization.
Further, 
\[\delta_{k-1}(G) \ge a -k \ge \frac12\left(n- C \sqrt{n}\right)-k \ge n/2- C\sqrt{n}, \]
since each $S \in \binom{A\cup B}{k-1}$ forms an edge with all vertices in either $A \setminus S$ or $B\setminus S$.
But $C$ was chosen depending only on $k$.
Therefore, we are done in this case.
\end{proof}
\begin{remark}
More precisely, our construction (or slight variations on it) can obtain the following based on the congruence class of $k$ modulo $4$:
\begin{itemize}
    \item If $k \equiv 0 \pmod{4}$, then the statement holds for infinitely many $n$ with the stronger conclusion $\delta_{k-1}(G) \ge \frac12 n - O(k)$.
    \item If $k \equiv 1,2 \pmod{4}$, then the statement holds for all sufficiently large $n$ with the stronger conclusion $\delta_{k-1}(G) \ge \frac12 n - O(k)$.
    \item If $k \equiv 3 \pmod{4}$, then the statement holds for all sufficiently large $n$.
\end{itemize}
It seems like the missing cases for $n$ and the weaker bounds on $\delta_{k-1}(G)$ are artifacts of the way we adapt the construction in~\cite{RRS:09} (adding matchings that are perfect on $A$), and it should be possible to prove the statement for all sufficiently large $n$ with the stronger conclusion $\delta_{k-1}(G) \ge \frac12n-O(k)$.
\end{remark}

\section{Concluding remarks} \label{sec:conclusion}
In this paper, we obtained a minimum degree version of Baranyai's theorem.
There are many intriguing related questions for future research.

While our paper determines the approximate $(k-1)$-degree threshold for a vertex-regular $k$-graph to admit a $1$-factorization, it would be interesting to determine the threshold for other degree types.
For $1\le \ell < k \le n$ with $k\mid n$, let $\delta(n,k,\ell)$ be the minimum $\delta$ such that for any vertex-regular $k$-graph $G$ on $n$ vertices with $\delta_\ell(G) \ge \delta$, $G$ admits a $1$-factorization.
Theorem~\ref{thm:main} together with Proposition~\ref{prop:lowerbound} yields $\limsup_{n\to\infty,\, k\mid n}\delta(n,k,k-1)/n = 1/2$.
\begin{problem}
    Given $1\le \ell <k$, determine $\limsup_{n\to\infty,\, k\mid n}\delta(n,k,\ell)/\binom{n}{k-\ell}$.
\end{problem}

It would also be interesting to determine the exact minimum $(k-1)$-degree threshold for $1$-factorizations of $k$-graphs, thus generalizing the result of Csaba, K\"uhn, Lo, Osthus, and Treglown~\cite{CBKLOT:16} to higher uniformities.
\begin{problem}
    Given $k \ge 3$ and $n$ sufficiently large compared to $k$, determine $\delta(n,k,k-1)$.
\end{problem}

The following conjecture would provide a Hamilton cycle analogue to our result.
A $k$-graph is a \defn{(tight) cycle} if it admits a cyclic ordering of its vertices such that its edges are the sets of $k$ consecutive vertices in this ordering, and a \defn{Hamilton cycle} in a $k$-graph $G$ is a cycle subgraph of $G$ that covers all vertices of $G$.
\begin{conj} \label{conj:hamilton}
    Given $k\ge 2$, $\delta>0$, $n$ large, and an $n$-vertex $r$-regular $k$-graph $G$ with $\delta_{k-1}(G) \ge (\frac12+\delta)n$, where $k \mid r$, $G$ admits a decomposition into Hamilton cycles.
\end{conj}
Because a Hamilton cycle in $G$ has a $1$-factorization when $k \mid n$, assuming Conjecture~\ref{conj:hamilton} holds, one can obtain a $1$-factorization of a vertex-regular Dirac $k$-graph $G$ with $k \mid |V(G)|$ by first removing a constant number of edge-disjoint perfect matchings so that the remaining vertex-degree is divisible by $k$ and then decomposing the remaining $k$-graph into Hamilton cycles, implying Theorem~\ref{thm:main}.
Conjecture~\ref{conj:hamilton} would also imply the following conjecture of Bailey and Stevens in the case that $n$ is large compared to $k$.
\begin{conj}[Bailey, Stevens~\cite{BS:10}]
    Let $k \ge 2$ and $n \ge k+2$ such that $n \mid \binom{n}{k}$.
    Then $K_n^{(k)}$ admits a decomposition into Hamilton cycles.
\end{conj}

One can also consider a $k$-partite analogue of Theorem~\ref{thm:main}.
Say that a $k$-graph $G$ is \defn{$k$-partite} with parts $V_1,\ldots,V_k$ if $V_1,\ldots,V_k$ form a partition of $V(G)$ and $|e\cap V_i|=1$ for all $e\in E(G)$ and $i\in [k]$.
For $1 \le \ell <k$, let $\delta'_\ell(G)$ denote the minimum of $d_G(S)$, taken over all $S\in \binom{V(G)}{\ell}$ with $|S\cap V_i| \le 1$ for all $i\in [k]$.
\begin{conj}
    Given $k\ge 2$, $\delta > 0$, $n$ large, and a regular $k$-partite $k$-graph $G$ with parts of size $n$, if $\delta'_{k-1}(G) \ge (\frac12+\delta)n$, then $G$ admits a $1$-factorization.
\end{conj}

Another direction for further research is enumerating the number of distinct $1$-factorizations a Dirac $k$-graph admits.
For $k=2$, this is known asymptotically: Ferber, Jain, and Sudakov~\cite{FJS:2020} proved that an $n$-vertex $d$-regular graph, where $n$ is even and $d \ge (\frac12+o(1))n$, admits $((1+o(1))\frac{d}{e^2})^{nd/2}$ distinct $1$-factorizations.

In a random $n$-vertex $k$-graph in which each edge appears independently with probability $p$, the expected number of perfect matchings is $\Phi\left(K_n^{(k)}\right) \cdot p^{n/k}$, where $\Phi\left(K_n^{(k)}\right) = \frac{n!}{(\frac{n}{k})!(k!)^{n/k}}$ denotes the number of perfect matchings in $K_n^{(k)}$.
By choosing perfect matchings one by one, we then heuristically expect an $r$-regular Dirac $k$-graph to admit close to $\left((1+o(1))\frac{r}{e^k}\right)^{nr/k}$ $1$-factorizations.
Thus, we have the following conjecture.
\begin{conj}
    Given $k\ge 2$, $\delta>0$, $n$ large with $k \mid n$, and an $r$-regular $k$-graph $G$ with $\delta_{k-1}(G) \ge (\frac12+\delta)n$, $G$ admits $\left((1+o(1))\frac{r}{e^k}\right)^{rn/k}$ distinct $1$-factorizations.
\end{conj}
The fact that $G$ admits at most $\left((1+o(1))\frac{r}{e^k}\right)^{rn/k}$ distinct $1$-factorizations follows from the entropy method.

\section*{Statement of AI use}
All mathematical ideas, arguments, results and proofs have been developed and written by the authors. This project was initiated in summer 2025 and a first draft with all essential arguments existed by spring 2026. In summer 2026, the main result was communicated at various conferences. When finishing the manuscript, we used ChatGPT-6 Astra only for local improvements of the~language and style of the manuscript, and for checking the mathematical correctness of all proofs.

\bibliographystyle{amsplain_initials_nobysame_nomr.bst}
\bibliography{References}

\end{document}